# On Problems Related to Primes: Some Ideas

Dhananjay P. Mehendale
Sir Parashurambhau College, Tilak Road, Pune-411030,
India

## Abstract

We present some new ideas on important problems related to primes. The topics of our discussion are: simple formulae for primes, twin primes, Sophie Germain primes, prime tuples less than or equal to a predefined number, and their infinitude; establishment of a kind of similarity between natural numbers and numbers that appear in an arithmetic progression, similar formulae for primes and the so called generalized twin primes in an arithmetic progression and their infinitude; generalization of Bertrand postulate and a Bertrand like postulate for twin primes; some elementary implications of a simple primality test, the use of Chinese remainder theorem in a possible proof of the Goldbach conjecture; Schinzel Sierpinski conjecture; and the Mersenne primes and composites, Fermat primes and their infinitude. Lastly we define other twin primes and provide a simple argument in support of their infinitude.

**1. Introduction:** One of the most remarkable theorems in number theory having a remarkably simple proof is due to **Euclid** [1]

**Theorem 1.1 (Euclid):** There are infinitely many primes.

**Proof:** There are infinitely many primes because odd multiples of finitely many odd primes do not cover all odd numbers. Every odd is of type $2n-1, n \in Z$, the set of positive integers. If there are finitely many primes $\{p_1(=2), p_2, \cdots, p_k\}$ then all odd numbers would be covered in the numbers of type $\{p_i(2n-1), 2 \leq i \leq k, n \in Z\}$ which is not true, since successive odd numbers are separated from each other by separation equal to 2 while these numbers are separated from one another by minimum separation equal to $2p_i, i > 1$, and if we add in the list of primes the uncovered numbers that are newly found and which have no smaller factors other than 1 or themselves then the same story repeats, i.e. we get each time new uncovered numbers, for example (as stated by Euclid) if one assumes that there are only finitely many primes say $p_1, p_2, \cdots, p_k$ then the number $\prod_{i=1}^{k} p_i + 1$ is either an uncovered number itself or it is a product of uncovered numbers other than $p_1, p_2, \cdots, p_k$, since it does not belong to the set $\{ p_i(2n-1), 2 \leq i \leq k, n \in Z \}$, hence etc. □

The problems related to primes which have received a wider attention in mathematical as well as nonmathematical circles are respectively **twin prime conjecture** and **Goldbach conjecture**. These two problems top the list of unsolved problems in number theory. In this paper we deal with these problems.

**2. Formulae for primes and twin primes $\leq x$ :** The number of primes and twin primes less than or equal to $x$, $\pi(x)$ and $T(x)$ respectively, can be expressed in terms of the following formulae: The following lemma 2.1 is a standard well-known result which directly follows from the sieve of **Eratosthenes** using the inclusion-exclusion principle [2].

**Lemma 2.1:** Let $S = \{ p_1, p_2, \cdots, p_k \}$ be set of first $k$ primes and $p_k^2 \leq x < p_{(k+1)}^2$, then

$$\pi(x) = x - \sum_{i=1}^{k} \left[ \frac{x}{p_i} \right] + \sum_{i<j}^{k} \left[ \frac{x}{p_i p_j} \right] - \sum_{i<j<l}^{k} \left[ \frac{x}{p_i p_j p_l} \right] + \cdots$$

$$\cdots + (-1)^k \left[ \frac{x}{p_1 p_2 \cdots p_k} \right] + (p_k - 1)$$

where [q] represents the integral part of q. □

If now we replace the square brackets around the numbers representing the integral part by ordinary round brackets and treat the numbers enclosed as real numbers and collect the terms (ignoring the last term) we get the following function:

$$\psi(x) = x \prod_{p_j \leq \sqrt{x}} (1 - \frac{1}{p_j})$$

It can be seen that as $x \to \infty$ $\psi(x)$, $\pi(x)$ change in such a way that **they keep approaching each other**, so that we can grant

$$\psi(x) \approx \pi(x) \qquad \to (2.1)$$

An alternative way to see the same thing is the following:
**Theorem 2.1:** $\psi(x) \sim \pi(x)$ in the asymptotic limit.

**Proof:** By prime number theorem $\pi(x) \sim \dfrac{x}{\log(x)}$. It is well-known that as $n \to \infty$, $\left(\left(\sum \dfrac{1}{n}\right) - \log(n)\right) \to \gamma (= 0.577\cdots)$, $\to$ Euler's constant.

i.e. $\sum \dfrac{1}{n} \sim \log(n)$. Now, $\sum \dfrac{1}{n} = \prod_p \dfrac{1}{(1 - \dfrac{1}{p})}$, $\to$ Euler's formula.

So, using the smallness of end correction and prime number theorem, we have

$$\psi(x) = x \prod_{p_j \leq \sqrt{x}} (1 - \dfrac{1}{p_j}) \sim \pi(x). \quad \square$$

The next lemma 2.2 can also be easily obtained, again, using inclusion-exclusion principle as described below:

**Lemma 2.2:** Let $S = \{p_1, p_2, \cdots, p_k\}$ be set of first $k$ primes and $p_k^2 \leq x < p_{(k+1)}^2$, then

$$T(x) = x + A + B + T(\sqrt{x}) \qquad \to (1)$$

where,

$$A = -\left[\dfrac{x}{2}\right] + \sum_{j=2}^{k} \sum_m \left(\left[\dfrac{x}{2p_j(m)}\right]\right) - \sum_{j<l\,m_1,m_2}^{k} \sum \left(\left[\dfrac{x}{2p_j(m_1)p_l(m_2)}\right]\right) +$$
$$\cdots + (-1)^k \sum_{m_i} \left[\dfrac{x}{2p_2(m_2)p_3(m_3)\cdots p_k(m_k)}\right]$$

$$B = -\sum_{j=2}^{k}\sum_{m}\left[\frac{x}{p_j(m)}\right] + \sum_{j<l}\sum_{m_i}\left(\left[\frac{x}{p_j(m_1)p_l(m_2)}\right]\right) - \cdots$$

$$+ (-1)^{(k-1)}\sum_{m_i}\left(\left[\frac{x}{p_2(m_2)\cdots p_k(m_k)}\right]\right)$$

where $m, m_i \in \{0,2\}$ for all $i$ and where $[q]$ represents the integral part of q. Also, note that $\left[\dfrac{x}{p_j(m)}\right]$ counts the numbers $\leq x$ which are congruent to $m \bmod(p_j)$ and $\left[\dfrac{x}{p_j(m_i)}\right]$ counts the numbers $\leq x$ which are congruent to $m_i \bmod(p_j)$.

**Proof:** Twin primes-pairs $\leq x$ are those primes-pairs $(p-2, p)$ for which $p \leq x$ is prime such that $(p-2)$ is also a prime. Thus, $(p-2, p)$ forms a twin-prime pair if and only if
(a) $p$ is **not** congruent to $0 \bmod(p_i)$, and (also, simultaneously)
(b) $p$ is **not** congruent to $2 \bmod(p_i)$.
where $p_i \in S$. In other words, the pair $(p-2, p)$ will form a twin-prime pair when $p \equiv \alpha_i \bmod(p_i)$, and $\alpha_i \neq 0, 2$ for all $i = 1, 2, \cdots, k$. Thus, in order to count the twin-prime pairs we need to count **those primes** $p$ which satisfy both the conditions (a) and (b). We proceed with this counting in a similar way as is done to find $\pi(x)$ in terms of the well known formula stated in the lemma 2.1 and obtain the above given formula for $T(x)$.

In counting something using inclusion-exclusion principle the term which is subtracted first (first term in $A$ and $B$) eliminates the undesired terms and the other terms are to bring about the corrections caused by under/over counting.

In the formula stated in equation (1), $A$ takes into account the terms that **do** involve the first prime $p_1(=2)$ in S, and $B$ takes into account the terms that **do not** involve this first prime. Since 2 =

0 in *GF* (2) (Galois field modulo 2) we have split the formula in terms of $A$ and $B$ for the sake of clarity and convenience.

(1) The first term in $A$ eliminates numbers $\leq x$ that are congruent to $0 \bmod(p_1(=2))$ while the first term in $B$ eliminates numbers $\leq x$ that are congruent to $0 \bmod(p_j)$ or $2 \bmod(p_j)$, $2 \leq j \leq k$.

(2) The numbers which are simultaneously congruent to $0 \bmod(p_1(=2))$ (first term in A) and also congruent to $0 \bmod(p_j)$ or $2 \bmod(p_j)$, $2 \leq j \leq k$ get subtracted twice. So, in order to bring about the correction, we have taken the second term with plus sign in $A$. Similarly, the numbers which are simultaneously congruent to $0 \bmod(p_j)$ or $2 \bmod(p_j)$, and $0 \bmod(p_l)$ or $2 \bmod(p_l)$, $2 \leq j < l \leq k$ get subtracted twice. So, in order to bring about the correction, we have taken the second term with plus sign in $B$.

(3) Continuing on similar lines with usual considerations of inclusion-exclusion principle one can easily see why there are other terms in $A$ and $B$ with respective specific signs. Also, note that Twin primes-pairs $\leq \sqrt{x}$ are not counted, in fact they get omitted during formation of A and B, so we add them separately. Hence, the lemma. □

**Example 2.1:** $T(20) = 20 + (-10+3+3) + (-6-6) = 4$. Note that the twin-prime pairs $\leq 20$ are $\{(3, 5), (5, 7), (11, 13), (17, 19)\}$.

For large $x$, $\left[\dfrac{x}{p_j(0)}\right] \approx \left[\dfrac{x}{p_j(2)}\right]$. So with this we can write the above formula for twin primes in a simplistic form as follows:

$$T(x) = x - \sum_{p_i \leq \sqrt{x}} 2\left[\frac{x}{p_i}\right] + \sum_{i<j} 2^2\left[\frac{x}{p_i p_j}\right] - \sum_{i<j<l} 2^3\left[\frac{x}{p_i p_j p_l}\right] + \ldots$$

$$\ldots + (-1)^k 2^k\left[\frac{x}{p_1 p_2 \cdots p_k}\right] + T(\sqrt{x})$$

where [q] represents the integral part of q.

Consider the following function $\omega(x)$ obtained from $T(x)$ by replacing the square brackets representing integer part of the number inside by ordinary brackets and ignoring the last term, thus:

$$\omega(x) = x - \sum_{p_i \leq \sqrt{x}} 2\left(\frac{x}{p_i}\right) + \sum_{i<j} 2^2\left(\frac{x}{p_i p_j}\right) - \sum_{i<j<l} 2^3\left(\frac{x}{p_i p_j p_l}\right) + \ldots$$

$$\ldots + (-1)^k 2^k \left(\frac{x}{p_1 p_2 \cdots p_k}\right),$$ which can be further written as

$$\omega(x) = \left(\frac{x}{2}\right) \prod_{2 < p \leq \sqrt{x}} (1 - \frac{2}{p})$$

It can be seen that as $x \to \infty$, $\omega(x), T(x)$ change in such a way that **they keep approaching each other**, so that we can grant

$$\omega(x) \approx T(x) \qquad \to (2.2)$$

**3. Probability Theory Arguments for Primes and Twin Primes:** We now proceed to discuss probability theory arguments that justify approximate equalities (2.1) and (2. 2).

**Theorem 3.1:** The number of primes $\leq x$, $\pi(x)$, are **approximately** equal to $\psi(x) = x \prod_{p_j \leq \sqrt{x}} (1 - \frac{1}{p_j})$.

**Proof:** A number $\alpha \leq x$ is prime if and only if $\alpha \equiv \beta_i \mod(p_i)$ and $\beta_i \neq 0$ for all primes $p_i \leq \sqrt{x}$. For a fixed $p_j$ the remainder $\beta_j$ can take one of the values $\{0, 1, \cdots, p_j - 1\}$. Therefore, for an $\alpha \leq x$ chosen at random, the probability that $\beta_j \neq 0$ will be $\left(\frac{p_j - 1}{p_j}\right)$. Clearly, for an $\alpha \leq x$ **chosen at random**, the probability that $\beta_i \neq 0$ for all $i = 1, 2, \cdots k$ such that $p_k^2 \leq x < p_{(k+1)}^2$ and $p_1, p_2, \cdots, p_k$ are all primes $\leq x$ will be $\prod_{i=1}^{k}\left(\frac{p_i - 1}{p_i}\right) = \prod_{i=1}^{k}\left(1 - \frac{1}{p_i}\right)$. Therefore, the number of primes less than or equal to $x$ will be

$$x\prod_{i=1}^{k}\left(1 - \frac{1}{p_i}\right) = x \prod_{p_i \leq \sqrt{x}}\left(1 - \frac{1}{p_i}\right) = \psi(x). \quad \square$$

**Theorem 3.2:** The number of twin primes $\leq x$, $T(x)$, are **approximately** equal to $\omega(x) = \frac{1}{2}x \prod_{2 < p_j \leq \sqrt{x}}(1 - \frac{2}{p_j})$.

**Proof:** A pair of numbers $(\alpha - 2, \alpha)$ forms a twin prime pair if both these numbers are prime. This can be equivalently stated as follows: a pair of numbers $(\alpha - 2, \alpha)$ forms a twin prime pair if $\alpha \equiv \beta_i \mod(p_i)$ and $\beta_i \neq 0, 2$ for all primes $p_i \leq \sqrt{x}$. For if some $\beta_i = 0$ then $\alpha$ will not be prime while if some $\beta_i = 2$ then $\alpha - 2$ will not be prime. In the case when $i = 1$, i.e. $p_i = p_1(= 2)$ then $\beta_i$ can take only two values, 0 and 1, and in this case for an $\alpha \leq x$ chosen at random, the probability that $\beta_1 \neq 0$ will be $\left(\frac{p_1 - 1}{p_1}\right) = \frac{1}{2}$. In all other cases for an $\alpha \leq x$ chosen at random, the probability that $\beta_j \neq 0, 2$ will be $\left(\frac{p_j - 2}{p_j}\right)$.

Clearly, for an $\alpha \leq x$ chosen at random, the probability that $\beta_i \neq 0, 2$ for all $i = 1, 2, \cdots k$ such that $p_k^2 \leq x < p_{(k+1)}^2$ and $p_1, p_2, \cdots, p_k$ are all primes $\leq x$ will be $\frac{1}{2}\prod_{i=2}^{k}\left(\frac{p_i - 2}{p_i}\right) = \frac{1}{2}\prod_{i=2}^{k}\left(1 - \frac{2}{p_i}\right)$. Therefore, the number of twin prime pairs less than or equal to $x$ will be

$$\frac{1}{2}x\prod_{i=2}^{k}\left(1 - \frac{2}{p_i}\right) = \frac{1}{2}x \prod_{2 < p_i \leq \sqrt{x}}\left(1 - \frac{2}{p_i}\right) = \omega(x). \quad \square$$

**Remark 3.1:** Theorems 3.1, 3.2 justify equalities (2.1), (2.2) respectively.

**4. Twin Primes:** If we show that $\omega(x) \to \infty$ as $x \to \infty$ then it follows from the observation $\omega(x) \approx T(x)$ that twin primes are infinite. We now proceed to accomplish this task by using methods of analytic number theory [3].

**Theorem 4.1:** Twin primes are infinite.

**Proof:** Let $U(x) = \prod_{2 < p \leq \sqrt{x}} (1 - \frac{2}{p})$.

Taking logs we have

$$\log U(x) = \sum \log(1 - \frac{2}{p})$$

In order to estimate this sum we consider the following power series:

$$-\log(1 - u) = u + \frac{u^2}{2} + \frac{u^3}{3} + \cdots$$

where $|u| < 1$. Using $u = \frac{2}{p}$ and taking the linear term to the left hand side we have

$$\phi(p) = -\log(1 - \frac{2}{p}) - \frac{2}{p}$$

It is clear that

$$0 < \phi(p) = (\frac{2^2}{2p^2} + \frac{2^3}{3p^3} + \cdots) < \frac{1}{2}(\frac{2^2}{p^2} + \frac{2^3}{p^3} + \cdots) = \frac{2}{p(p-2)}$$

The above inequality shows that the infinite series $M = \sum \phi(p)$ converges since it is dominated by series $\sum_n \left(\frac{2}{n(n-2)}\right)$.

Thus,

$$0 < M - \sum_{2 < p \leq \sqrt{x}} \phi(p) = \sum_{p > \sqrt{x}} \phi(p) \leq \sum_{n \geq \sqrt{x}} \left(\frac{2}{n(n-2)}\right).$$

Hence,

$$\sum_{2 < p \leq \sqrt{x}} \phi(p) = M + O\left(\frac{1}{x}\right).$$

Thus,

$$-\log U(x) = \sum_{2 < p \leq \sqrt{x}} \left(\frac{2}{p}\right) + M + O\left(\frac{1}{x}\right).$$

Now, one can easily see that

$$\sum_{2 < p \leq \sqrt{x}} \left(\frac{2}{p}\right) = 2 \sum_{2 < p \leq \sqrt{x}} \left(\frac{1}{p}\right) = 2 \log \log(\sqrt{x}) = 2A + O\left(\frac{1}{\log(\sqrt{x})}\right)$$

where $A = 1 - \log\log(2) + \int_2^\infty \frac{R(t)}{t \log^2(t)} dt$, and $R(x) \sim O(1)$.

Therefore,

$$-\log U(x) = 2\log\left(\frac{1}{2}\log(x)\right) + 2A + O\left(\frac{1}{\log(x)}\right) + M + O\left(\frac{1}{x}\right)$$

So,

$$\log U(x) = -\log\left(\frac{1}{2}\log(x)\right)^2 - 2A - M - O\left(\frac{1}{\log(x)}\right)$$

Therefore,

$$U(x) = \frac{C}{\left(\frac{1}{2}\log(x)\right)^2} = \frac{4C}{(\log(x))^2}$$

Hence,

$$\omega(x) = \frac{1}{2} x U(x) \sim \frac{2Cx}{(\log(x))^2},$$

where, $C = \exp(-M - 2A)$

Now, using the observation $\omega(x) \approx T(x)$ we have $T(x) \approx \dfrac{2Cx}{(\log(x))^2}$, therefore, twin primes are **infinite**! $\square$

**Remark 4.1:** Let $S = \{p_1, p_2, \cdots, p_k\}$ be set of first $k$ primes and $p_k^2 \leq x < p_{(k+1)}^2$. The Sophie Germain-prime pairs $\leq x$ are those primes-pairs $(p, 2p+1)$ for which $p \leq x$ is prime such that $(2p+1)$ is also a prime. Thus, $(p, 2p+1)$ forms a Sophie Germain-prime pair if and only if
(a) $p$ is **not** congruent to $0 \bmod(p_i)$, and (also, simultaneously)
(b) $p$ is **not** congruent to $\beta \bmod(p_i)$ such that $2\beta + 1 \equiv 0 (\bmod(p_i))$.
where $p_i \in S$.

It is important to note that modulo a prime $p_i \in S$ under consideration there exists a **unique** remainder $\beta$ apart from the remainder 0 which is disallowed to achieve primality of the pair $(p, 2p+1)$. So, here also, we can proceed on similar lines and can develop **similar theorems** for the **Sophie Germain primes** like the theorems that are developed above for **twin primes**.

We now define and study a special kind of **Dirichlet Function**, $\xi_T(s)$.

$$\xi_T(s) = \sum_{n=1}^{\infty} \left( \dfrac{2^{\phi(n)}}{n^s} \right),$$ where $n$ has the following prime factorization:

$n = p_1^{m_1} p_2^{m_2} \cdots p_r^{m_r}$ and $\phi(n) = (m_1 + m_2 + \cdots + m_r)$.

We now record some simple **observations:**

**Lemma 4.1:** $\xi_T(1) = \sum_{n=1}^{\infty} \left( \dfrac{2^{\phi(n)}}{n} \right)$ diverges.

**Proof:** Since $\log(n) \sim \sum_{n=1}^{\infty} \dfrac{1}{n} < \xi_T(1)$ and $\log(n)$ is divergent as $n \to \infty$. $\square$

**Lemma 4.2:** $\psi(n) = 2^{\phi(n)}$ is totally multiplicative.

**Proof:** Clear. Since if $k = mn$ then $\psi(k) = \psi(m)\psi(n)$. □

**Lemma 4.3:** $\xi_T(1) = \prod_p (1 - \frac{2}{p})^{(-1)}$, where the product is taken over all primes $p$.

**Proof:** When $p$ is prime its prime factorization is $p^1$ and so $\psi(p) = 2^1 = 2$. Now,

$$\prod_p (1 - \frac{2}{p})^{(-1)} = (1 + \frac{2}{p_1} + \frac{2}{p_1^2} + \cdots)(1 + \frac{2}{p_2} + \frac{2}{p_2^2} + \cdots) \cdots$$

$$= \sum \frac{2^{\alpha_1 + \alpha_2 + \cdots \alpha_k}}{p_1^{\alpha_1} p_2^{\alpha_2} \cdots p_k^{\alpha_k}} = \sum \frac{2^{\phi(n)}}{n} = \xi_T(1). □$$

**Lemma 4.4:** The Euler product for $\xi_T(s)$ is

$$\prod_p (1 - \frac{\psi(p)}{p^s})^{(-1)} = \prod_p (1 - \frac{2}{p^s})^{(-1)}$$

where the product is taken over all primes $p$.

**Proof:** Expanding the Euler product we have

$$1 + 2\sum_p \frac{1}{p^s} + 2^2 \sum_{p_1 \leq p_2} \frac{1}{(p_1 p_2)^s} + \cdots + 2^k \sum_{p_1 \leq \cdots \leq p_k} \frac{1}{(p_1 \cdots p_k)^s} \text{ and it is}$$

clearly equal to $\sum \frac{2^{\phi(n)}}{n^s}$. □

**Lemma 4.5:** $T(x) = \frac{x}{2} \xi_T^*(1)$, where $\xi_T^*(1)$ denotes the Euler product for $\xi_T(1)$ taken over all primes up to $\sqrt{x}$.

**Proof:** Clear from theorem 3.2. □

Let $0 < \sigma < 1$.

**Lemma 4.6:** $\log \xi_T(1 + \sigma) = \sum_p \frac{2}{p^{(1+\sigma)}} + O(1)$.

**Proof:** For $0 < x < 1$,

$$-\log(1-x) = \sum_{n=1}^{\infty}\left(\frac{x^n}{n}\right).$$ Since $\sigma > 0$ we have

$$\xi_T(1+\sigma) = \prod_p (1 - \frac{2}{p^{(1+\sigma)}})^{(-1)} > 1, \text{ and}$$

$$\log \xi_T(1+\sigma) = \log \prod_p (1 - \frac{2}{p^{(1+\sigma)}})^{(-1)}$$

$$= -\sum_p \log(1 - \frac{2}{p^{(1+\sigma)}})$$

$$= \sum_p \sum_{n=1}^{\infty} \frac{2^n}{np^{n(1+\sigma)}}$$

$$= \sum_p \frac{2}{p^{(1+\sigma)}} + \sum_p \sum_{n=2}^{\infty} (\frac{2^n}{np^{n(1+\sigma)}}), \text{ therefore,}$$

$$\log \xi_T(1+\sigma) = \sum_p \frac{2}{p^{(1+\sigma)}} + O(1). \ \square$$

**Lemma 4.7:** $\log \xi_T(1+\sigma) = \log(\frac{1}{\sigma}) + O(1)$

**Proof:** From $0 < \sigma < 1$,

$$1 < \frac{1}{\sigma} = \int_1^{\infty} \frac{1}{x^{(1+\sigma)}} dx < \xi_T(1+\sigma) < 1 + \int_1^{\infty} \frac{1}{x^{(1+\sigma)}} dx = \frac{1}{\sigma} + 1, \text{ and so,}$$

$$0 < \log(\frac{1}{\sigma}) < \log \xi_T(1+\sigma) < \log(\frac{1}{\sigma}) + 1, \text{ therefore,}$$

$$\log \xi_T(1+\sigma) = \log(\frac{1}{\sigma}) + O(1). \ \square$$

The connection between $T(x)$ and $\xi_T^*(1)$ can now be utilized to establish the infinitude of twin-prime numbers. Equating the right hand sides of in the lemmas 6 and 7, we have

$$\log(\frac{1}{\sigma}) = \sum_p \frac{2}{p^{(1+\sigma)}} + O(1).$$

Now, if we assume that twin primes are finite then as $\sigma \to 0$ the right hand side of the above equation will remain bounded but the left hand side goes to infinity which is absurd. Hence, twin primes must be infinite.

**5. Prime k–Tuples Conjecture:** A $(k-1)$–tuple $(b_1, b_2, \cdots, b_{(k-1)})$ is called **admissible** when

(a) $b_1 < b_2 < \cdots < b_{(k-1)}$, and

(b) for every prime $q \leq k$ the residue classes $\{0 \bmod q, b_1 \bmod q, \cdots, b_{(k-1)} \bmod q\}$ is properly contained in the set of all residue classes modulo $q$, (see page 200 of the extremely good book, [5]). Among the admissible $(k-1)$–tuples those with minimal $b_{(k-1)}$ are said to be **tight**. Thus for $k \leq 4$ the tight admissible $(k-1)$–tuples are (2), (2, 6), (4, 6), (2, 6, 8).

While deriving a formula for twin primes we had noticed and used the fact that a pair of numbers $(n, n+2)$ forms a twin prime pair if and only if $n + 2 \equiv \alpha_i \bmod(p_i)$, $i = 1, 2, \cdots l$ and $p_l^2 < n < p_{(l+1)}^2$ such that $\alpha_i \neq 0, 2$. Similarly, if $k \geq 2$ and $(b_1, b_2, \cdots, b_{(k-1)})$ is an admissible $(k-1)$–tuple of positive integers, then numbers $p, p + b_1, p + b_2, \cdots, p + b_{(k-1)}$ are simultaneously prime if and only if $p + b_k \equiv \alpha_i \bmod(p_i)$, $i = 1, 2, \cdots l$ and $p_l^2 < (p + b_k) < p_{(l+1)}^2$ such that $\alpha_i \neq 0, b_1, b_2, \cdots, b_{(k-1)} (\bmod(p_i))$.

We now describe the so called **appropriate values** of remainders and **forbidden values** of remainders. The **appropriate values** of remainders are those which when occur do not cause discarding of the tuple under consideration for its nonprimality. The **forbidden values** of remainders are those which when occur force us to discard the tuple in the counting for its nonprimality. Let $u$ be the cardinality of the set of **forbidden values** of remainders. When the cardinality of the set of remainders, $(= p_j)$, corresponding to a prime $p_j$ is smaller than the cardinality of the set of admissible values then the appropriate values of remainders $(\bmod(p_j))$ are not $k$ in number but are **less**. As an **example** consider the set of admissible values mentioned hereinbefore, namely, (2, 6, 8). The smallest set of primes corresponding to this set is {11, 13, 17, 19}. For this set of admissible values the **forbidden values** of remainders are {0, 2, 6, 8}. Now when $p_j = 2$, the set of forbidden values of remainders ($\mod(p_j = 2)$) is {0}, so, from the set of remainders {0, 1} ($\mod(p_j = 2)$) the value {0} is forbidden (and so $u = 1$). Note that the set of the appropriate values of remainders ($\mod(p_j = 2)$) is {1}. When $p_j = 3$, the set of forbidden values of remainders ($\mod(p_j = 3)$) is {0, 2}, so, from the remainder set of remainders {0, 1, 2} the values {0, 2} are forbidden so $u = 2$. Note that the set appropriate values of remainders ($\mod(p_j = 3)$) is {1}. When $p_j = 5$, the set of values of forbidden values of remainders ($\mod(p_j = 5)$) is {0, 1, 2, 3}, so, from the set of remainders {0, 1, 2, 3, 4} the values {0, 1, 2, 3} are forbidden and so $u = 4$. In this case, the set appropriate values of remainders ($\mod(p_j = 5)$) is {4}. When $p_j = 7$, the set of forbidden values of remainders ($\mod(p_j = 7)$) is {0, 1, 2, 6}, so, from the set of remainders {0, 1, 2, 3, 4, 5, 6} the values {0, 1, 2, 6} are forbidden and so $u = 4$. Note that the set appropriate values of remainders ($\mod(p_j = 7)$) is {3, 4, 5} for this case. It is important to note for this example that when $p_j =$ **7 or higher**, proceeding on similar lines, we can see that $u = 4 = k$.

**Lemma 5.1:** Let $S = \{p_1, p_2, \cdots, p_l\}$ be set of first $l$ primes and $p_l^2 \leq x < p_{(l+1)}^2$, and let then the number of prime $k \geq 2$ and $(b_1, b_2, \cdots, b_{(k-1)})$ is an admissible $(k-1)$–tuple of positive integers, then the number of $(k)$–tuples of primes $(p, p+b_1, p+b_2, \cdots, p+b_{(k-1)}) \leq x$, $T_k(x)$, can be given by the following formula:

$$T_k(x) = x - \sum_{j=1}^{l} \sum_{m_i} \left[\frac{x}{p_j(m_i)}\right] + \sum_{j<k}^{l} \sum_{m_i} \left(\left[\frac{x}{p_j(m_i)p_k(m_j)}\right]\right) - \cdots$$

$$+ (-1)^l \sum_{m_i} \left(\left[\frac{x}{p_1(m_1)p_2(m_2)\cdots p_l(m_l)}\right]\right) + T_k(\sqrt{x}).$$

where $m_i \in \{0, b_1, b_2, \cdots, b_{(k-1)}\}$ for all $i$ (to be **precise**, $m_i \in$ the set of forbidden values of remainders) and where [q] represents the integral part of q. Also, $\left[\dfrac{x}{p_j(m_i)}\right]$ counts the numbers $\leq x$ which are congruent to $m_i \bmod(p_j)$. Note that prime k-tuples $\leq \sqrt{x}$ are counted in $T_k(\sqrt{x})$.

**Proof:** Proceeding on similar lines as is done for lemma 2.1 one can easily obtain the result. □

For large $x$, $\left[\dfrac{x}{p_j(m_i)}\right] \approx \left[\dfrac{x}{p_j(m_k)}\right]$. So with this we can write the above formula for the count of prime $k$-tuples in a simplistic form as follows:

$$T_k(x) = x - \sum_{p_i \leq \sqrt{x}} u\left[\dfrac{x}{p_i}\right] + \sum_{i<j} u^2 \left[\dfrac{x}{p_i p_j}\right] - \sum_{i<j<k} u^3 \left[\dfrac{x}{p_i p_j p_k}\right] + \ldots$$

$$\ldots + (-1)^l u^l \left[\dfrac{x}{p_1 p_2 \cdots p_l}\right] + T_k(\sqrt{x}).$$

where [q] represents the integral part of q, and
(i) $u = k$ when $p_i > b_{(k-1)}$.
(ii) $u = $ the cardinality of the set of forbidden values of remainders corresponding to the $p_i \leq b_{(k-1)}$ under consideration.

Now, consider the following function $\omega_k(x)$ obtained from $T_k(x)$ by replacing the square brackets representing integer part of the number inside by ordinary brackets, thus:

$$\omega_k(x) = x - \sum_{p_i \leq \sqrt{x}} u\left(\dfrac{x}{p_i}\right) + \sum_{i<j} u^2 \left(\dfrac{x}{p_i p_j}\right) - \sum_{i<j<k} u^3 \left(\dfrac{x}{p_i p_j p_k}\right) + \ldots$$

$$\ldots + (-1)^l u^l \left(\dfrac{x}{p_1 p_2 \cdots p_l}\right) + T_k(\sqrt{x}),$$ which can be further written as

$$\omega_k(x) = \prod_{p_i \leq \sqrt{x}} (1 - \frac{u}{p_i}) \qquad \rightarrow (5.1)$$

where [q] represents the integral part of q, and
(i) $u = k$ when $p_i > b_{(k-1)}$.
(ii) $u =$ the cardinality of the set of forbidden values of remainders corresponding to the $p_i \leq b_{(k-1)}$ under consideration.

It can be seen that as $x \rightarrow \infty$, $\omega_k(x), T_k(x)$ change in such a way that **they keep approaching each other**, so that we can grant

$$\omega_k(x) \approx T_k(x) \qquad \rightarrow (5.2)$$

As is done in section 3, theorem 3.2, for the case of twin primes, one can proceed **exactly on similar lines** with developing **probability theory arguments** to count the prime k-tuples $\leq \sqrt{x}$ that will lead to a formula justify approximate equality (5.2), namely, The **number of prime k-tuples** $\leq x$ are **approximately** equal to $\omega_k(x)$ and where

$$\omega_k(x) = \prod_{p_i \leq \sqrt{x}} (1 - \frac{u}{p_i})$$

where [q] represents the integral part of q, and
(i) $u = k$ when $p_i > b_{(k-1)}$.
(ii) $u =$ the cardinality of the set of forbidden values of remainders corresponding to the $p_i \leq b_{(k-1)}$ under consideration.

**Theorem 5.1(Prime k-Tuples Conjecture, (page 201, [5])):** If $k \geq 2$ and $(b_1, b_2, \cdots, b_{(k-1)})$ is an admissible $(k-1)$–tuple of positive integers, there exist infinitely many primes $p$ such that $p, p+b_1, p+b_2, \cdots, p+b_{(k-1)}$ are simultaneously prime.

**Proof:** As is done in section 4, if we show that $\omega_k(x) \rightarrow \infty$ as $x \rightarrow \infty$ then it follows from the observation $\omega_k(x) \approx T_k(x)$ that the above mentioned prime k-tuples in the statement of the theorem are infinite. We can proceed to accomplish this task by using methods of analytic number theory [3], on similar lines as is done in section 4 in the case of twin primes, by starting this time with equation (5.1), and show assuming the (justified) validity (by probability theory arguments) of equation (5.2) that

$$T_k(x) \approx \frac{\Omega.x}{(\log(x))^k},$$ where $\Omega$ is a constant. Hence etc. □

**Remark 5.1:** The **Hardy-Littlewood conjecture** states that if $\pi(x)$ denotes the number of primes $\leq x$ then $\pi(x+y) \leq \pi(x) + \pi(y)$, $x, y \geq 2$. Consider the number line:

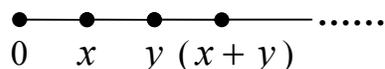

0   x   y  (x+y)

Clearly, $\pi(x+y) = \pi(x) + \pi(x \leftrightarrow (x+y))$  ….(1)
where $\pi(x \leftrightarrow (x+y))$ represents the count of primes between the interval $(x, (x+y)]$. Thus, the conjecture reduces to showing that for any fixed $y$, we have $\pi(x \leftrightarrow (x+y)) \leq \pi(y)$. Thus, the conjecture states that if we choose any finite interval, $(0, y)$, starting with zero on the number line and shift it anywhere, say $(x, x+y)$, on the number line and count the primes contained in the shifted interval then they are always less than or equal to the primes in the initial interval. This conjecture looks quite reasonable but the beautiful theorem of **Hensley and Richards** forbids its validity when the prime *k*-tuples conjecture is true (see page 202 of the very remarkable book, [5]).

**6. Generalized Bertrand Postulate:** Bertrand postulate assures the existence of a prime between *n* and 2*n* when n≥1. We now proceed to state a **generalization:**

**Theorem 6.1 (Generalized Bertrand Postulate):** For every $\alpha \in (1, 2]$, i.e. $1 < \alpha \leq 2$, there exists a positive integer $n_0$ such that for all $n \geq n_0$ there exists a prime *p* between *n* and $\alpha n$, i.e. $n < p \leq [\alpha n]$, where [q] represents the integral part of q.

**Proof:** By the prime number theorem, $\pi(x) \sim \dfrac{x}{\log(x)}$

Now, by proceeding as is done in [4] it can be seen that

$$\left(\frac{\alpha x}{\log(\alpha x)}\right) - \left(\frac{x}{\log(x)}\right) \sim \left(\frac{\alpha x}{(\log(x) + \log(\alpha))} - \frac{x}{\log(x)}\right) \sim$$

$$\left(\frac{x}{\log(x)} - \frac{\alpha x \log(\alpha)}{\log^2(x)}\right)$$

as $x \to \infty$. Thus, the primes in the interval between $x$ and $\alpha x$ are of the same order as those below $x$, hence etc. □

**Problem 6.1:** Discover some nontrivial implications of this generalization.

With the estimation of the of twin primes, $T(x) \approx \dfrac{2Cx}{(\log(x))^2}$ we can have a similar result like Bertrand Postulate for twin primes as well as follows:

**Theorem 6.2 (Bertrand Postulate for Twin Primes):** If $x$ is a real number, and $x \geq 7$, then there exists at least one twin-prime pair between $(x, 2x)$.

**Proof:** Considering $T(2x) - T(x)$, where $T(x) \approx \dfrac{2Cx}{(\log(x))^2}$, and using the similar considerations used to settle usual Bertrand postulate [4], one can easily get the result. □

**Theorem 6.3 (Generalized Bertrand Postulate for Twin Primes):** For every $\alpha \in (1, 2]$, i.e. $1 < \alpha \leq 2$, there exists a positive integer $n_0$ such that for all $n \geq n_0$ there exists a twin prime pair $(p, p+2)$ between $n$ and $\alpha n$, i.e. $n < p \leq [\alpha n]$, where [q] represents the integral part of q.

**Proof:** Follows using [4] by proceeding on similar lines as in theorem 4.1 and using $T(x) \approx \dfrac{2Cx}{(\log(x))^2}$. □

**7. Remainders $\mod(p)$ of numbers in Arithmetic Progression:** Note that $k + mp \equiv k \mod(p)$ for $0 \leq k \leq (p-1)$ and for all $m = \{0, 1, 2, \ldots\}$. Thus, remainders $\mod(p)$ of integers in the set of all integers, Z, form a **periodic sequence** $(0,1,\cdots,(p-1),0,1,\cdots(p-1),0,1,\cdots)$. Let $(a,b) = 1$ i.e. $a$, $b$ are co prime, and let $\{a+kb / k = 0, 1, 2, \ldots, t, \ldots\}$ be the set of numbers in the arithmetic progression. It is clear to see the interesting fact that remainders $\mod(p)$ of the numbers in the arithmetic progression form a **periodic sequence** with period $p$ made up of some **permutation of** numbers $\{0,1,\cdots,(p-1)\}$ when b $\neq mp$ for some integer $m$. Thus, let $(a + kb) \equiv \alpha_k \mod(p)$, $k = 0, 1, 2, \ldots, t, \ldots$ then it is easy to check that the sequence of remainders $\mod(p)$ of the numbers in this arithmetic progression for some fixed prime $p$ is either

(a) $(\alpha_0, \alpha_1, \cdots, \alpha_{(p-1)}, \alpha_0, \alpha_1, \cdots)$ i.e. a periodic sequence of some permutation of numbers $\{0,1,\cdots, p-1\}$ i.e. it repeats after every $p$ numbers, or

(b) $(\alpha_0, \alpha_1, \cdots, \alpha_{(p-1)}, \alpha_0, \alpha_1, \cdots) = (\beta, \beta, \cdots, \beta, \beta, \beta, \cdots)$ where $\beta$ is constant $(\neq 0)$, $b = mp$ for some integer $m$, and $a \equiv \beta \mod(p)$.

**8. Probability Theory Arguments for Primes and Twin Primes in an Arithmetic Progression:** As seen above, as far as the remainders modulo primes are concerned the numbers in any arithmetic progression are **same** as natural numbers and either the remainders only get **permuted** or they are a **nonzero constant** modulo a prime. The only essential change is the reduction of quantity of numbers. Thus, up to $x$ there are less numbers in the arithmetic progression. For example in $\{a+kb / k = 0, 1, 2, \ldots, t, \ldots\}$, $(a,b) = 1$, there are approximately $\left(\dfrac{x-a}{b}\right)$ numbers up to $x$. So, with **similar proofs** given for theorem 3.1, theorem 3.2 we will have the following:

**Theorem 8.1:** Let S denote the set of numbers in an arithmetic progression, $S = \{a+kb / k = 0, 1, 2, \ldots, t, \ldots\}$, and $(a,b) = 1$. The number of primes $\leq x$ in S is **approximately** equal to

$$\phi(x) = \left(\dfrac{x-a}{b}\right) \prod_{p_j \leq \sqrt{x}} (1 - \dfrac{1}{p_j}). \quad \square$$

**Theorem 8.2:** Let S denote the set of numbers in an arithmetic progression, $S = \{a+kb / k = 0, 1, 2, \ldots, t, \ldots\}$, and $(a,b) = 1$. The number of generalized twin prime pairs, of type
$((a + kb),(a + (k + 1)b))$ when all numbers in the arithmetic progression or odd, and of type $((a + kb),(a + (k + 2)b))$ when numbers in the arithmetic progression or successively odd and even, $\leq x$ in S is **approximately** equal to

$$\theta(x) = \frac{1}{2}\left(\frac{x-a}{b}\right) \prod_{2 < p_j \leq \sqrt{x}} (1 - \frac{2}{p_j}) \quad \square$$

By proceeding on similar lines as is done for theorem 2.1, and theorem 4.1 we can have the following theorems:

**Theorem 8.3:** Let S denote the set of numbers in an arithmetic progression, $S = \{a+kb / k = 0, 1, 2, \ldots, t, \ldots\}$, and $(a,b) = 1$. The number of primes $\leq x$ in S, say $\pi_b(x)$, is **approximately** equal to

$$\pi_b(x) \approx \left(\frac{((x-a)/b)}{\log((x-a)/b)}\right) \quad \square$$

**Theorem 8.4:** Let S denote the set of numbers in an arithmetic progression, $S = \{a+kb / k = 0, 1, 2, \ldots, t, \ldots\}$, and $(a,b) = 1$. The number of generalized twin prime pairs, of type $((a + kb),(a + (k + 1)b))$ when all numbers in the arithmetic progression or odd and of type $((a + kb),(a + (k + 2)b))$ when numbers in the arithmetic progression or successively odd and even, $\leq x$ in S is **approximately** equal to

$$T_b(x) \approx \left(\frac{2C((x-a)/b)}{(\log((x-a)/b))^2}\right),$$

where $C$ is the constant similar to one obtained in the theorem 3.1. $\square$

**Definition 8.1:** If $p$ does not divide $a$ and the smallest $k$ satisfying $a^k \equiv 1 \bmod(p)$ is $(p-1)$ then $a$ is called the primitive root $\bmod(p)$.

We take an opportunity to state the following **conjecture** at this juncture:

**Conjecture 8.1:** For every integer $Q$ other than $1, -1$, and or perfect square there exist infinitely many primes $p$ for which Q is primitive root modulo $p$ in every arithmetic progression $a + kb / k = 0, 1, 2, \ldots, t, \ldots$ such that $(a, b) = 1$, and

$$N_b^Q(x) \approx A(Q)\left(\frac{(x-a)/b}{\log((x-a)/b)}\right),$$

where $N_b^Q(x)$ denote the number of such primes less than or equal to $x$ for which Q is primitive root modulo $p$, and $A(Q)$ is certain constant depending on $Q$. □

**9. Goldbach Conjecture:** It asserts that every even integer greater than 2 is the sum of two primes. Stated in a letter to Leonard Euler by Christian Goldbach in 1842, this is still an enduring unsolved problem. There are many other conjectures about prime sums given in [6].

We begin our discussion with some elementary observations. If a number $m$ is not prime and $m = pq$ then clearly either $p \leq \sqrt{m}$ or $q \leq \sqrt{m}$. Further, if $p \leq \sqrt{m}$ say, and suppose $p$ is not prime then by fundamental theorem of arithmetic $p$ can be factored uniquely as the product of prime powers with obviously all the primes in that unique factorization strictly less than $p$. Thus, we have the following well-known

**Theorem 9.1:** A number $m$ is prime if and only if it is not divisible by any prime number $\leq \sqrt{m}$. □

Let $S = \{p_1, p_2, \cdots, p_k\}$ be the set of first $k$ primes and
$$p_k^2 \leq N < p_{(k+1)}^2$$
Consider any disjoint partitioning of set of primes in $S$ into two parts, like
$A = \{p_{11}, p_{12}, \cdots p_{1m}\}$ and $B = \{p_{21}, p_{22}, \cdots p_{2n}\}$ such that
$A \cap B = \phi$ and $A \cup B = S$.

**Preposition 9.1:** Let $\alpha = (p_{11}^{\mu_1} \cdot p_{12}^{\mu_2} \cdots p_{1m}^{\mu_m}) \pm (p_{21}^{\nu_1} \cdot p_{22}^{\nu_2} \cdots p_{2n}^{\nu_n})$ then $\alpha$ is prime if $2 < |\alpha| < N$.

**Proof:** Let $P = p_{11}^{\mu_1} \cdot p_{12}^{\mu_2} \cdots p_{1m}^{\mu_m}$ and $Q = p_{21}^{\nu_1} \cdot p_{22}^{\nu_2} \cdots p_{2n}^{\nu_n}$, then P is not divisible by prime in the set $B$ though Q is and Q is not divisible by a prime in the set A though P is, therefore P ± Q is not divisible by a prime in S. Also, $2 < |\alpha| < N$, hence by theorem 1.1, $\alpha$ is prime. □

**Preposition 9.2:** Let $A = \{p_{11}, p_{12}, \cdots, p_{1m}\}$, $B = \{p_{21}, p_{22}, \cdots, p_{2n}\}$, and $C = \{p_j\}$ such that $A, B, C$ are mutually disjoint sets, and $A \cup B \cup C = S$. Let $p_1 (= 2)$ belongs to A. Let $\alpha = \prod_{l=1}^{m} p_{1l}^{\delta_l} - p_j$. Now,

let $R = \prod_{l=1}^{m} p_{1l}^{\delta_l}$, then

(a) If $2 < |\alpha| < N$ and

(b) If $R \mod (p_{2k}) \neq p_j \mod (p_{2k})$ for all $1 \leq k \leq n$ then $\alpha$ is prime.

**Proof:** Since R is divisible by primes in the set A while $p_j$ is not, and $p_j$ is divisible by $p_j$ but R is not so $\alpha$ is not divisible by primes in $A \cup C$. Also, from (b) $\alpha$ is not divisible by any prime in B, therefore $\alpha$ is not divisible by any prime in S and (a) holds therefore by theorem 1.1 $\alpha$ is prime. □

Note that R is even and $R = \alpha + p_j$, i.e. sum of two primes.

**Preposition 9.3:** Let $N = 2n$ and let $T = \{p_1, p_2, \cdots, p_k\}$ be set of first $k$ primes and are $\leq n$. Let $N \equiv \beta_j \mod(p_j)$, $p_j \in T$. Let $C = \{p_{(k+1)}, \cdots, p_l\}$ be the set of next primes after $p_k$ and $p_s \in C$ satisfies $n \leq p_s < 2n = N$, then $N - p_s$ is prime if $p_s \equiv \alpha_j \mod(p_j)$ and $\alpha_j \neq \beta_j$ for all primes $p_j \in T$.

**Proof:** It is clear from the given data that $(N - p_s) \equiv (\beta_j - \alpha_j) \mod(p_j)$ and so for primality of $N - p_s$ we should have $\alpha_j \neq \beta_j$ for all $p_j \in T$. □

Note that when such a prime $p_s \in C$ exists then N is expressible as a sum of two primes.

**Preposition 9.4:** There are infinitely many even numbers that can be expressed as sum of two primes.

**Proof:** Let $p$ be some fixed odd prime. Since primes are infinite in number so obviously there are infinitely many odd primes $q > p$ and clearly $p + q$ is an even number. □

We now proceed to propose a **fresh new** approach to deal with Goldbach problem and finally show that this new strategy is worthwhile for affirmative settlement of the problem. We begin with the statement of the following **very useful** theorem:

**Theorem 9.2 (Chinese Remainder Theorem):** Let $m_1, m_2, \cdots, m_l$ denote $l$ positive integers which are relatively primes in pairs, and let $a_1, a_2, \cdots, a_l$ denote any $l$ integers. Then the following congruence system

$$x \equiv a_1 \mod(m_1)$$
$$x \equiv a_2 \mod(m_2)$$
$$\vdots$$
$$x \equiv a_l \mod(m_l)$$

has common solutions. If $x_0$ is one such solution, then an integer $x$ is another solution if and only if $x = x_0 + km$ for some integer $k$ and here $m = \prod_{j=1}^{l} m_j$. □

Let $2n$ be an arbitrary positive and even integer. In order to settle Goldbach conjecture in the affirmative we have to show the existence of two prime numbers $p$, $q < 2n$ such that $2n = p + q$. Let $p_1(= 2), p_2(= 3), p_3, \cdots, p_k, p_{(k+1)}$ be the successive primes such that $p_k^2 < 2n < p_{(k+1)}^2$. Let $2n \equiv \beta_i \mod(p_i), i = 1, 2, \cdots, k$. Note that since $2n$ is even so clearly $\beta_1 = 0$. The **main idea** in this approach is to split each remainder $\beta_i$ into **two nonzero parts in all possible ways:** $\beta_i = \eta_i^j + \delta_i^j$, and $1 \leq \eta_i^j, \delta_i^j \leq (p_i - 1)$. We then determine all possible number pairs ($p, q$) using the above mentioned Chinese remainder theorem such that $p \equiv \eta_i^j \mod(p_i)$ and $q \equiv \delta_i^j \mod(p_i)$. We call such $p, q$ numbers the **suitable candidates**. They are **complements** of each other in the sense that $\beta_i = \eta_i^j + \delta_i^j$. If we show that among the suitable candidates there exists at least one number $p$ (or $q$) such that $p < 2n$ (and so prime due to theorem 9.1) then Goldbach conjecture follows. We show the existence of such prime $p$ by showing that the **span** (separation between the largest and the smallest number among the suitable candidates) of the suitable candidates, $\prod_{i=1}^{k}(p_i - u_i)$ in count, where $u_i = 1$ when $\beta_i = 0$, and $u_i = 2$ when $\beta_i \neq 0$, is greater than $\prod_{i=1}^{k} p_i - 2n$.

For the sake of **illustration** we begin with an example:
**Example 9.1:** Let $2n = 100$. So $\sqrt{2n} = 10$. The primes less than 10 are respectively 2, 3, 5, and 7. Cleary, $100 \equiv 0 \bmod(2)$, $100 \equiv 1 \bmod(3)$, $100 \equiv 0 \bmod(5)$, and $100 \equiv 2 \bmod(7)$. Now we note down all the possibilities that exists for $\eta_i^j$:

(1) $\eta_1^j = \{1\}$
(2) $\eta_2^j = \{2\}$
(3) $\eta_3^j = \{1,2,3,4\}$
(4) $\eta_4^j = \{1,3,4,5,6\}$

Note that for any choice of $\eta_i^j$ given above the corresponding choice for $\delta_i^j$ that get fixed by the requirement, namely, $\beta_i = \eta_i^j + \delta_i^j$, is suitable in the sense that all these $\delta_i^j$ are nonzero as is needed.

To settle Goldbach conjecture for the even number equal to 100 we need numbers (**at least one**) less than 100 which can be expressed simultaneously in the forms
$2k_1 + 1, 3k_2 + 2, 5k_3 + \{1 or 2 or 3 or 4\}, 7k_4 + \{1 or 3 or 4 or 5 or 6\}$ for some positive integers $k_1, k_2, k_3, k_4$. For example, note that 11 has the desired representations i.e. 2×5+1, 3×3+2, 5×2+1, 7×1+4. etc. so, 11 is the suitable prime with suitable prime complement 89 so that 100 = 11+89.

Note that there are in all 1×1×4×5 = 20 (product of the cardinalities of sets $\eta_i^j$, $i = \{1,2,3,4\}$) possibilities for $p$ and one can see that the numbers obtained from these choices by applying **Chinese remainder theorem** are {11, 17, 29, 41, 47, 53, 59, 71, 83, 89, 101, 113, 131, 137, 143, 167, 173, 179, 197, 209}, when written in increasing order. Thus, there are in all 10 choices less than $2n(=100)$ for $p$ and since $p + q = q + p$ therefore there are 5 ways to express 100 as sum of two primes, namely, 11+89, 17+83, 29+71, 41+59, 47+53. Here we have split the remainders $\beta_i$ into two nonzero parts.

But instead if we allow the splitting such that some one $\eta_i^j = 0$ and $\delta_i^j \neq 0$ then it leads to additional expressions for 100 as sum of two primes, namely, 3+97. Here we allow a prime among primes 2, 3, 5, 7 as one prime when its complement is also prime.

Now, note that $\eta_1^j = \{1\}$ and $\eta_2^j = \{2\}$, so they have only one choice. So, all the possible 20 numbers are separated from each other by at least 2×3 = 6. When $\eta_1^j, \eta_2^j, \eta_3^j$ are fixed and only $\eta_4^j$ is allowed to take all possible values mentioned above then it is easy to see that the resulting 5 numbers are separated from each other by at least 2×3×5 = 30. Similarly, when $\eta_1^j, \eta_2^j, \eta_4^j$ are fixed and only $\eta_3^j$ is allowed to take all possible values mentioned above then it is easy to see that the resulting 4 numbers are separated from each other by at least 2×3×7 = 42. Thus, one can easily see that the mutual separations among the allowed 20 numbers constructed above are either multiples of 6 or 30 or 42. The following are the conclusions about the 20 numbers that are possible as candidates for $p$:

(1) All these 20 numbers have relative separation equal to a nonzero integral multiple of 6.
(2) There are four groups of numbers (since we can keep fixed the first three remainders in 4 ways) such that each group contains 5 numbers and the numbers in each such a group have a relative separation equal to a nonzero integral multiple of 30.
(3) There are five groups of numbers (since we can keep fixed the first, second and fourth remainder in 5 ways) such that each group contains 4 numbers and the numbers in each such a group have a relative separation equal to a nonzero integral multiple of 42.

Now, suppose the maximum number among the twenty numbers is of the order of $M = \prod_{i=1}^{4} p_i$.

(a) Thus, among the four groups described in (2) there is a group having largest element $M$. It is easy to see that the largest element in other groups will lower down at least by 2×3 = 6. Now, since there are four such groups of type described in (2) there will be a group of elements having largest element equal to $M - 24$. Since the elements in this group are separated by 30 we will have an element among the twenty element of the order $M - 24 - 150 = 36$. Thus, the span between the largest number, say $M = \prod_{i=1}^{4} p_i$ and smallest number which is of order $(M - 24 - 150)$ will be 174, while $M - 2n = \prod_{i=1}^{4} p_i - 2n$ is of the order 110.

(b) Thus, among the five groups described in (3) there is a group having largest element $M$. It is easy to see that the largest element in other groups will lower down at least by $2 \times 3 = 6$. Now, since there are five such groups of type described in (3) there will be a group of elements having largest element equal to $M - 30$. Since the elements in this group are separated by 42 we will have an element among the twenty element of the order $M - 24 - 168 = 12$. Thus, the **span** between the largest number, say $M = \prod_{i=1}^{4} p_i$ and smallest number, which is of order $(M - 24 - 168)$, will be 192, while $M - 2n = \prod_{i=1}^{4} p_i - 2n$ is of the order 110.

Note that the actual largest number in the list of the solutions obtained above using Chinese remainder theorem is 209 and therefore the smallest number in the list that we get is $209 - 192 = 11$.

Let us now proceed with some **useful lemmas:**

**Lemma 9.1:** Let $2n$ be an arbitrary positive and even integer and let $p_1(=2), p_2(=3), p_3, \cdots, p_k, p_{(k+1)}$ be the successive primes such that $p_k^2 < 2n < p_{(k+1)}^2$ and let $2n \equiv \beta_i \mod(p_i), i = 1, 2, \cdots, k$. Then $\beta_i$ can be expressed as sum of two nonzero integers respectively in exactly $p_i - 1$ ways when $\beta_i = 0$ and in exactly $p_i - 2$ ways when $\beta_i \neq 0$.

**Proof:** (1) If $\beta_i = 0$ then $\beta_i$ can be expressed as sum of two nonzero integers as follows: $\beta_i = 0 = 1 + (p_i - 1)$
$$= 2 + (p_i - 2)$$
$$= 3 + (p_i - 3)$$
$$\vdots$$
$$= (p_i - 1) + 1.$$

Thus, in all exactly $p_i - 1$ ways.

(2) If $\beta_i \neq 0$ then $\beta_i$ can be expressed as sum of two nonzero integers as follows: Let $\beta_i = j$, $1 \leq j \leq (p_i - 1)$ then

$$\begin{aligned} \beta_i = j &= 1 + (j-1) \\ &= 2 + (j-2) \\ &\vdots \\ &= (j-1) + 1 \\ &= (j+1) + (p_i - 1) \\ &= (j+2) + (p_i - 2) \\ &\vdots \\ &= (p_i - 1) + (j+1). \end{aligned}$$

Thus, in all exactly $p_i - 2$ ways. Hence the lemma. □

**Lemma 9.2:** Let $2n \equiv \beta_i \mod(p_i), i = 1, 2, \cdots, k$. The total number of possibilities for number $p$ such that $p \equiv \eta_i^j \mod(p_i)$, $q \equiv \delta_i^j \mod(p_i)$ and $\beta_i = \eta_i^j + \delta_i^j$, where $1 \leq \eta_i^j, \delta_i^j \leq (p_i - 1)$ are N in number and $N = \prod_{i=1}^{k}(p_i - u_i)$ where $u_i = 1$ when $\beta_i = 0$, and $u_i = 2$ when $\beta_i \neq 0$.

**Proof:** Since all the choices for $\eta_i^j$ for all $i$ can be made independently of each other therefore the result follows from lemma 9.1. □

**Lemma 9.3:** If we fix the choice of some $j$ remainders among the suitable remainders corresponding to some $j$ primes and allow to vary the other remainders over the range of suitable remainders corresponding to other $(k - j)$ primes and generate all numbers using Chinese remainder theorem then these numbers will be separated from each other by minimum separation equal to $\prod_l p_l$ where the product is taken over the selected $j$ primes such that the choice for their remainders is fixed.

**Proof:** It follows from the periodic occurrence of remainders. A fixed remainder $\eta_i^j$ for prime $p_i$ reappears after $p_i$ numbers, similarly, the simultaneous occurrence of fixed remainders $\eta_i^j$ and $\eta_l^k$ for primes $p_i$ and $p_l$ together, takes place after $p_i \cdot p_l$ numbers. Hence etc. □

**Lemma 9.4:** The span of suitable candidates in the group of $(p_k - u_k)$ numbers obtained using Chinese remainder theorem by choosing fixed remainders $\eta_i^j$, $i = 1, 2, \cdots, (k-1)$ and allowing to vary $\eta_k^j$ over the allowed $(p_k - u_k)$ remainders is at least equal to $\prod_{i=1}^{k-1} p_i (p_k - u_k)$, i.e. if the largest element in the group of these $(p_k - u_k)$ numbers is of the order $M = \prod_{i=1}^{k} p_i$ then the smallest element is at least of the order

$$M - \prod_{i=1}^{k-1} p_i (p_k - u_k).$$

**Proof:** Since the minimum separation between any two numbers in this group of numbers is $\prod_{i=1}^{k-1} p_i$ and there are in all $(p_k - u_k)$ numbers. Hence etc. □

**Lemma 9.5:** Suppose the largest number in the group of $(p_k - u_k)$ elements obtained by fixing the choice of allowed remainders $\eta_i^j$, $i = 1, 2, \cdots, (k-1)$ and allowing to vary $\eta_k^j$ over the allowed $(p_k - u_k)$ remainders is of order $M$ (i.e. $\approx M = \prod_{i=1}^{k} p_i$) then we can have a group of $(p_k - u_k)$ numbers having largest element of the order

$$M - \sum_{j=1}^{(k-3)} \prod_{i=1}^{(k-j-1)} p_i (p_{(k-j)} - u_{(k-j)}).$$

**Proof:** Suppose we have a group of $(p_k - u_k)$ elements obtained by fixing the choice for $\eta_i^j$, $i = 1, 2, \cdots, (k-1)$ and allowing to vary $\eta_k^j$, having largest element of the order M. In this case if we allow to vary $\eta_{(k-1)}^j$ along with $\eta_k^j$ then the largest element can be lowered to $M - \prod_{i=1}^{k-2} p_i (p_{(k-1)} - u_{(k-1)})$. Proceeding on similar lines if we also allow varying $\eta_{(k-2)}^j$ then the largest element can further be lowered to

$$M - \prod_{i=1}^{k-2} p_i (p_{(k-1)} - u_{(k-1)}) - \prod_{i=1}^{k-3} p_i (p_{(k-1)} - u_{(k-1)}).$$

Since for primes $p_1 (= 2), p_2 (= 3)$ there is unique choice for remainders so successively allowing to vary $\eta_i^j$, $i = 3, 4, \cdots, (k-1)$, we can ultimately lower the largest element to L where

$$L = M - \sum_{j=1}^{(k-3)} \prod_{i=1}^{(k-j-1)} p_i (p_{(k-j)} - u_{(k-j)}). \quad \square$$

**Lemma 9.6:** The span of suitable candidates for an even number $2n$ ranges from $M = \prod_{i=1}^{k} p_i$ to $(L - \prod_{i=1}^{k-1} p_i (p_k - u_k))$, where L is as defined in the above lemma 9.5.

**Proof:** Clear, since there exists a group of $(p_k - u_k)$ elements separated from each other at least by $\prod_{i=1}^{k-1} p_i$ and having largest element of the order L. $\square$

**Remark 9.1:** Thus, the smallest suitable candidate is of the order of

$$M - \sum_{j=0}^{(k-3)} \prod_{i=1}^{(k-j-1)} p_i (p_{(k-j)} - u_{(k-j)}).$$

**Lemma 9.7:** If we assume that all $u_i = 2$ then the smallest element among the suitable candidates will be

$$M - \sum_{j=0}^{(k-3)} \prod_{i=1}^{(k-j-1)} p_i (p_{(k-j)} - 2) = \prod_{i=1}^{k-1} p_i + \cdots + \prod_{i=1}^{k-3} p_i + \ldots + 2 p_1 \cdot p_2$$

**Proof:** Substituting $M = \prod_{i=1}^{k} p_i$ the result is clear by simple evaluation. $\square$

**Lemma 9.8:** If we assume that all $u_i = 1$ then the smallest element among the suitable candidates will be

$$M - \sum_{j=0}^{(k-3)} \prod_{i=1}^{(k-j-1)} p_i (p_{(k-j)} - 1) = p_1 \cdot p_2 = 6.$$

**Proof:** Substituting $M = \prod_{i=1}^{k} p_i$ the result is clear by simple evaluation. □

**Remark 9.2:** Note that when $u_k = 2$ then in this case one remainder is disallowed because it's choosing leads to a zero remainder for one of the two solutions ($p$ or $q$) through Chinese remainder theorem and thus implying its nonprime nature. But it is **important to note** that this number causing zero remainder in general do not have a place at the lower end of the set of numbers separated from each other by the product of primes whose remainder choices are fixed, the probability for such an occurrence is equal to $\left(\dfrac{1}{(p_k - 1)}\right)$ and so the probability that the omitted number producing zero remainder does not have a place at the lower end will be equal to $\left(\dfrac{(p_k - u_k)}{(p_k - 1)}\right)$, and so the smallest element could be of the order $M - \sum_{j=1}^{(k-3)} \prod_{i=1}^{(k-j-1)} p_i (p_{(k-j)} - 1)$.

In the above **example 9.1**, if we fix the remainders $\eta_i^j$, $i = 1,2,3$ as {1, 2, 2} and allow to vary the remainder $\eta_i^j$ for $i = 4$ over the entire range of remainders form 1 to 6, including the forbidden remainder (= 2), then we get the following numbers using the Chinese remainder theorem:

| Remainder sequence mod(2,3,5,7) | The corresponding number using Chinese remainder theorem |
|---|---|
| (1, 2, 2, 1) | 197 |
| **(1, 2, 2, 2)** | **107** |
| (1, 2, 2, 3) | 17 |
| (1, 2, 2, 4) | 137 |
| (1, 2, 2, 5) | 47 |
| (1, 2, 2, 6) | 167 |

Thus, the number (= 107) corresponding to the forbidden remainder (= 2) modulo(7) **has not appeared at the lower end** of the span.

**Theorem 9.3:** The span of suitable candidates for an even number $2n$ is **greater** than $(\prod_{i=1}^{k} p_i - 2n)$ where $p_1(=2), p_2(=3), p_3, \cdots, p_k, p_{(k+1)}$ are the successive primes such that $p_k^2 < 2n < p_{(k+1)}^2$, i.e. there exists at least one suitable candidate, $p$ say, such that $p < 2n$.

**Proof:** With the consideration of high probability of nonoccurrence of the omitted choices producing zero remainder at the lower end we can replace almost all $u_i$ in the above span by $1$ and we get the span of the order $\approx \sum_{j=1}^{(k-3)(k-j-1)} \prod_{i=1} p_i (p_{(k-j)} - 1) + \prod_{i=1}^{k-1} p_i (p_k - 1)$.

and this span is clearly greater than $\prod_{i=1}^{k} p_i - 2n$. (Note that $2n \approx p_k^2$).

Hence etc. □

**Corollary (Goldbach Conjecture):** Every even integer greater than 2 is the sum of two primes.

**Proof:** The existence of $p < 2n$ as per above theorem implies the existence of its complement $q < 2n$ such that $2n = p + q$.
□

The above given proof is **existential** type. To obtain the actual representation for an even number as sum of two primes when we land up at a number $> 2n$ in the list of suitable candidates the following strategy is **sometimes** useful:
  (1) Express the suitable candidate we have arrived at which is $> 2n$, $t$ say, as $t = n + r$ (here, $r > n$).
  (2) Check whether we can factorize $r = sr'$ such that
      $s^2 \equiv 1 \bmod(p_i)$ for all $i = 1, 2, \cdots, k$. (It is clear that such $s$ will satisfy the congruence relations $s \equiv \omega_i \bmod(p_i)$, where $\omega_i = 1$ or $\omega_i = p_i - 1$.)
  (3) If yes, then check whether $r' < n$, if yes then $(n - r', n + r')$ will form the desired prime pair whose sum is $2n$.

**Example 9.2:** In example 9.1 we get number 137 among the suitable candidates. Now, $137 = 50 + 87$. Now, $87 = 29 \times 3$ and $29^2 \equiv 1 \mod(p_i)$ for all $i = 1, 2, \cdots, 4$. So, $(50-3, 50+3)$ forms the desired prime pair such that the sum is 100.

Twin prime conjecture can be viewed through **similar considerations** like Goldbach conjecture. Affirmative settlement of twin prime conjecture is equivalent to showing the existence of infinitely many odd numbers $n$ such that if $n \equiv \alpha_i \mod(p_i)$, $i = 1, 2, \cdots l$ and $p_l^2 < n < p_{(l+1)}^2$ such that for all $i$, either $\alpha_i \neq 0, 2$ or $\alpha_i \neq 0, (p_i - 2)$. The validity of the first case implies that $(n-2, n)$ is a twin prime pair while the validity of the second case implies that $(n, n+2)$ is a twin prime pair.

**Example 9.3:** Consider all numbers modulo 2, 3, 5 up to $2 \times 3 \times 5 = 30$ numbers. The sufficient condition (as per theorem 9.1) among these numbers with $\alpha_1 = \{1\}$, $\alpha_2 = \{1\}$, $\alpha_3 = \{1, 3, 4\}$ to form twin prime pairs $(n-2, n)$ is $n < 7^2 = 49$. Using Chinese remainder theorem we see that these numbers "$n$" are $\{13, 19, 31, 43\}$ when arranged in increasing order, and so we get four twin prime pairs of the form $(n-2, n)$, namely, $\{(11, 13), (17, 19), (29, 31), (41, 43)\}$.

**Example 9.4:** Consider all numbers modulo 2, 3, 5, 7 up to $2 \times 3 \times 5 \times 7 = 210$ The sufficient condition (as per theorem 9.1) among these numbers with $\alpha_1 = \{1\}$, $\alpha_2 = \{2\}$, $\alpha_3 = \{1, 2, 4\}$, $\alpha_4 = \{1, 2, 3, 4, 6\}$ to form twin prime pairs of the type $(n, n+2)$ is $n < (11)^2 = 121$. Using Chinese remainder theorem we see that these numbers "$n$" are $\{11, 17, 29, 41, 59, 71, 101, 107, 137, 149, 167, 179, 191, 197, 209\}$ when arranged in increasing order, and so **with certainty** the following eight are twin prime pairs of the form $(n, n+2)$, namely, $\{(11, 13), (17, 19), (29, 31), (41, 43), (59, 61), (71, 73), (101, 103), (107, 109)\}$.

One can see that considering all numbers "$n$" modulo primes $p_1, p_2, \cdots, p_k$ from 1 to $\prod_{i=1}^{k} p_i$ such that if $n \equiv \alpha_i \mod(p_i)$, $i = 1, 2, \cdots k$ and for all $i$, $\alpha_i \neq 0, 2$ as solutions as per Chinese remainder theorem, with pre-specified remainders to get the twin prime pairs of type $(n-2, n)$, (then applying **similar considerations** as is done for the settlement of Goldbach conjecture) one can see that the count of twin prime pairs of type $(n-2, n)$ such that $p_k^2 < n < p_{(k+1)}^2$ goes on **increasing** with increase in $k$. Thus, one is assured to have endless supply of twin prime pairs.

**Alternative Treatment for Goldbach Conjecture:** We now proceed to show how this famous Goldbach conjecture actually follows as a consequence of Chinese remainder theorem stated above in quite transparent way. For this let us start with the following equivalent:

**Theorem 9.4 (Equivalent of Goldbach Conjecture):** For every positive integer greater than or equal to 5 there exists two primes equidistant from it. In other words, let $n \geq 5$ be the positive integer then there exists distance $d$, $0 < d < n$, such that $p = n - d$ and $q = n + d$, where $p, q$ are prime numbers.

The equivalence of the above statement with Goldbach conjecture is straightforward. If the above statement (equivalent) is valid then it implies that $n - p = q - n = d$ which in turn implies that $2n = p + q$, which is Goldbach Conjecture.

Let $p_1 (= 2), p_2 (= 3), p_3, \cdots, p_k, p_{(k+1)}$ be the successive primes such that $p_k^2 < 2n < p_{(k+1)}^2$. Let $n \equiv \alpha_i \mod(p_i), i = 1, 2, \cdots, k$ and further each $\alpha_i$ satisfies the inequality $0 \leq \alpha_i \leq (p_i - 1)$. Now, in order to find out the desired primes $p, q$ mentioned above which are equidistant from $n \geq 5$ and at distance $d$, $0 < d < n$, we need to have following congruence relations, namely, $n - d \equiv \mu_i \mod(p_i)$, and $n + d \equiv \nu_i \mod(p_i)$, such that $\mu_i > 0, \nu_i > 0$ for all $i = 1, 2, \cdots, k$. We now proceed to record the **formulae of distances** for the case of each prime and for each remainder modulo that prime. Modulo each prime, $p_i$ and for each possible value of $\alpha_i$ for that prime $p_i$ for each such that at those distances there will be positive entries (positive values for $\mu_i, \nu_i$).

Let us record below few cases concretely:

1) **Mod(2) case:** Let $\alpha_1 = 0$ then for numbers at distance $d = 2k_1 + 1$, on both left and right side of the number $n$ on the number line, we will have positive remainders modulo 2, where, $k_1 = 0,1,2,\cdots$ Similarly, Let $\alpha_1 = 1$ then at distance $d = 2k_1$, on both left and right side of the number $n$ on the number line, we will have positive entries, where, $k_1 = 0,1,2,\cdots$

2) **Mod(3) case:** Let $\alpha_2 = 0$ then for numbers at distance $d = 3k_2 + 1, 3k_2 + 2$, on both left and right side of the number $n$ on the number line, we will have positive remainders modulo 3, where, $k_2 = 0,1,2,\cdots$ Similarly, Let $\alpha_2 = 1,2$ then for numbers at distance $d = 3(k_2 + 1)$, on both left and right side of the number $n$ on the number line, we will have positive remainders modulo 3, where, $k_2 = 0,1,2,\cdots$

3) **Mod(5) case:** Let $\alpha_3 = 0$ then for numbers at distance $d = 5k_3 + 1, 5k_3 + 2, 5k_3 + 3, 5k_3 + 4$, on both left and right side of the number $n$ on the number line, we will have positive remainders modulo 5, where, $k_3 = 0,1,2,\cdots$ Similarly: Let $\alpha_3 = 1$ then for numbers at distance $d = 5k_3 + 2, 5k_3 + 3, 5(k_3 + 1)$, on both left and right side of the number $n$ on the number line, we will have positive remainders modulo 5, where, $k_3 = 0,1,2,\cdots$, Let $\alpha_3 = 2$ then for numbers at distance $d = 5k_3 + 1, 5k_3 + 4, 5(k_3 + 1)$, on both left and right side of the number $n$ on the number line, we will have positive remainders modulo 5, where, $k_3 = 0,1,2,\cdots$, Let $\alpha_3 = 3$ then for numbers at distance $d = 5k_3 + 1, 5k_3 + 4, 5(k_3 + 1)$, on both left and right side of the number $n$ on the number line, we will have positive remainders modulo 5, where, $k_3 = 0,1,2,\cdots$, Let $\alpha_3 = 4$ then for numbers at distance $d = 5k_3 + 2, 5k_3 + 3, 5(k_3 + 1)$, on both left and right side of the number $n$ on the number line, we will have positive remainders modulo 5, where, $k_3 = 0,1,2,\cdots$.

One can continue in this way and determine the formulae for distance $d$, $0 < d < n$, for all primes $p_1(=2), p_2(=3), p_3, \cdots, p_k, p_{(k+1)}$ such that we will have positive remainders modulo corresponding prime under consideration for numbers at the distance $d$ satisfying these formulae, on both left and right side of the number *n* on the number line, and which (we aim to) should ultimately lead to the following congruence relations, namely, $n - d \equiv \mu_i \mod(p_i)$, and $n + d \equiv v_i \mod(p_i)$, such that $\mu_i > 0, v_i > 0$ for all $i = 1, 2, \cdots, k$, which settles Goldbach conjecture. Thus, in order settle Goldbach conjecture we need to find distance $d$, $0 < d < n$, which takes the form $d = p_i k_i + r_i$, for all primes $p_1(=2), p_2(=3), p_3, \cdots, p_k, p_{(k+1)}$ and $r_i \in \{0, 1, \cdots, p_i - 1\}$ such that we will have positive remainders modulo corresponding prime under consideration for numbers at the distance $d$ satisfying these formulae, on both sides, left and right side of the number *n* on the number line.

**Proof of theorem 9.4:** Now, given $n \geq 5$ we can determine uniquely the remainders $\alpha_i$ satisfying the inequality $0 \leq \alpha_i \leq (p_i - 1)$ such that $n \equiv \alpha_i \mod(p_i), i = 1, 2, \cdots, k$. Using these values of $\alpha_i$ we can find formulae (expressions) for distance $d$, as is concretely done above for Mod(2), Mod(3), Mod(5) cases, in the form $d = p_i k_i + r_i$, for all primes $p_1(=2), p_2(=3), p_3, \cdots, p_k, p_{(k+1)}$ and where we will get $r_i \in \{0, 1, \cdots, p_i - 1\}$ such that we will have positive remainders modulo corresponding primes under consideration for numbers at the distance $d$ satisfying these formulae, on both sides, left as well as right side of the number *n* on the number line. But what does these formulae for distances imply? These formulae for distances $d = p_i k_i + r_i$ equivalently imply that we in possession of following **congruence system** (and are after finding a positive number $d$, $0 < d < n$, satisfying this congruence system):

$$d \equiv r_1 \mod(p_1)$$
$$d \equiv r_2 \mod(p_2)$$
$$\vdots$$
$$d \equiv r_k \mod(p_k)$$

Now, clearly, this system of congruence **has a solution by Chinese Remainder Theorem**. Further using this solution, $d$, we can find

$p = n - d$ and $q = n + d$, where $p, q$ are prime numbers. And thus clearly we have the desired expression $2n = p + q$, for given even number, $2n$, as required for settling Goldbach Conjecture.

□

**Example 9.5:** Let $2n = 34$, therefore we have $n = 17$. Now, it is clear to see that clearly, $17 \equiv 1 \bmod(2)$, $17 \equiv 2 \bmod(3)$, and $17 \equiv 2 \bmod(5)$. Therefore for Mod(2) case: $\alpha_1 = 1$ and so $d = 2k_1$. Similarly, for Mod(3) case: $\alpha_2 = 2$ and so $d = 3k_2$. Similarly, for Mod(5) case: $\alpha_3 = 2$ and so $d = 5k_3 + 1$. Therefore, we have
$$d \equiv 0 \bmod(2)$$
$$d \equiv 0 \bmod(3)$$
$$d \equiv 1 \bmod(5)$$
so as a solution of this congruence system we get $d = 6$. Therefore, it implies that $n - p = q - n = d$ equivalent to $17 - 11 = 23 - 17 = 6$ which in turn from relation $2n = p + q$ implies that $34 = 11 + 23$.

**Example 9.6:** Let $2n = 100$, therefore we have $n = 50$. Now, it is clear to see that clearly, $50 \equiv 0 \bmod(2)$, $50 \equiv 2 \bmod(3)$, $50 \equiv 0 \bmod(5)$, and $50 \equiv 1 \bmod(7)$. Therefore for Mod(2) case: $\alpha_1 = 0$ and so $d = 2k_1 + 1$. Similarly, for Mod(3) case: $\alpha_2 = 2$ and so $d = 3k_2$. Similarly, for Mod(5) case: $\alpha_3 = 0$ and so $d = 5k_3 + \{1,2,3,4\}$. Similarly, for Mod(7) case: $\alpha_4 = 1$ and so $d = 7k_4 + \{2,3,4,5\}$. Therefore, we have the following congruence system:
$$d \equiv 1 \bmod(2)$$
$$d \equiv 0 \bmod(3)$$
$$d \equiv \{1,2,3,4\} \bmod(5)$$
and
$$d \equiv \{0,2,3,4,5\} \bmod(7)$$

For this congruence system we can easily find solutions, viz: $d = \{3, 9, 21, 33, 39\}$ and this leads to expression for 100 as follows: $100 = \{53+47, 59+41, 71+29, 83+17, 89+11\}$.

**One more Argument for Goldbach Conjecture:** We now offer an elementary counting argument that will imply the validity of Goldbach conjecture in most convincing manner. This argument has been inspired by the observation that if we will note down all primes up to certain positive integer $x$, say $p_1, p_2, \cdots, p_k$, where $k = \left[\dfrac{x}{\ln(x)}\right]$ denotes the integer part of $\left(\dfrac{x}{\ln(x)}\right)$, and consider all possible expressions like $p_1 + p_1$, $p_1 + p_2$, $p_1 + p_3$, ... , $p_1 + p_k$, $p_2 + p_2$, $p_2 + p_3$, ..., $p_k + p_k$, clearly these expressions represent some even numbers $\leq 2x$. Now, total number of even numbers $\leq 2x$ are $\sim x$, and total number of above expressions representing some (not all) even numbers $\leq 2x$ are $\sim \dfrac{k(k+1)}{2}$.

Clearly, $\dfrac{k(k+1)}{2} >> x$ as $x \to \infty$ and therefore not only many even numbers should have representation as sum of two primes but also those even numbers will have many such distinct representations.

The so called counting argument that we develop below in support of Goldbach conjecture goes as follows:

I. We count the number of distinct representations for chosen fixed positive and even integer $x$ as sum of two odd numbers bigger than one. **Let $A$ denotes this count**. . In these representations all cases like, e.g. both the odd numbers in the representation are primes, etc. are also included.
II. We count the number of representations for chosen fixed positive and even integer $x$ as sum of two odd numbers in which either both or at least one odd number in the representation is not prime. **Let $B$ denotes this count**
III. It is straightforward to see that if we subtract count obtained in II above from the count obtained in I then that results in the count for chosen even integer $x$ as sum of two prime numbers.

Let $x$ be the given even number. To find the count in I we need to find number of odd numbers starting from 3 up to $x/2$.

If x /2 is odd then $A = \left[\dfrac{x}{4}\right]$.

If x/2 is even then $A = \left[\dfrac{x}{4}\right] - 1$.

Note that square bracket [ ] denotes the integer part of.

**Chinese Remainder Theorem and Twin Prime Conjecture:**

We now proceed with direct convincing proof that establishes the desired infinitude of twin prime pairs as a consequence of Chinese Remainder Theorem.

**Theorem 9.5 (Twin Prime Conjecture):** Twin prime pairs are infinite.

**Proof:** Consider any sufficiently large number $N$. Suppose $p_1(=2), p_2(=3), p_3, \cdots, p_k, p_{(k+1)}$ are all successive primes such that $p_k^2 \leq N < p_{(k+1)}^2$. We now choose $x$ such that $x = 2k_1 + 1$, also $x = 3k_2 + 1$, $x = 5k_3 + \{1,3,4\}$, $x = 7k_4 + \{1,3,4,5,6\}$, $x = 11k_5 + \{1,3,4,5,6,7,8,9,10\}$ ,……, $x = p_k k_k + \{1,3,4,\cdots, p_k - 1\}$ and $k_i, i = 1,2,\cdots$ are suitable positive integers. These integers and corresponding remainders are so chosen such that $x \leq N$. Equivalently, we are finding such $x$ which is solution of the following congruence system:

$$x \equiv 1 \mod(2)$$
$$x \equiv 1 \mod(3)$$
$$x \equiv \{1,3,4\} \mod(5)$$
$$x \equiv \{1,3,4,5,6\} \mod(7)$$
$$x \equiv \{1,3,4,5,6,7,8,9,10\} \mod(11)$$
$$\vdots$$
$$x \equiv \{1,3,4,\cdots, p_k - 1\} \mod(p_k)$$

The existence of such $x \leq N$ ensures its prime nature and for such $x$ not only it will be prime but also it is easy to check that $(x-2)$ will also be a prime number! Thus, whichever endlessly large $N$ we choose we always can find a twin prime pair near it! Thus, twin primes are infinite!!

**Remark 9.1:** Following on the lines of theorem 9.5 given above, it is now a straightforward exercise to establish the infinitude of prime pairs separated by $2n$, where $n$ is any positive integer.

**Example 9.7:** Let $x = 7$. So, $x \equiv 1 \mod(2)$, $x \equiv 1 \mod(3)$. Clearly, $(x-2) = 5$ is prime. Let $x = 31$. So, $x \equiv 1 \mod(2)$, $x \equiv 1 \mod(3)$, $x \equiv 1 \mod(5)$. Clearly, $(x-2) = 29$ is prime. Let $x = 43$. So, $x \equiv 1 \mod(2)$, $x \equiv 1 \mod(3)$, $x \equiv 3 \mod(5)$. Clearly, $(x-2) = 41$ is prime. Let $x = 103$. So, $x \equiv 1 \mod(2)$, $x \equiv 1 \mod(3)$, $x \equiv 3 \mod(5)$, $x \equiv 5 \mod(7)$, Clearly, $(x-2) = 101$ is prime.

**10. Schinzel-Sierpinski Conjecture:** This well-known conjecture [7] asserts that every positive rational number $x$ can be represented as a quotient of shifted primes, i.e. there exist primes $p$ and $q$ such that

$$x = \left(\frac{p+1}{q+1}\right).$$

Let $x = \left(\frac{m}{n}\right)$.

Then one has to show the existence of a positive integer $k$ such that the pair $(p, q) = (2mk - 1, 2nk - 1)$, is a prime pair. Does there exist a **possible candidate** $k$ so that the pair $(p, q) = (2mk - 1, 2nk - 1)$ will turn out a prime pair? We proceed with an **algorithm** for searching such a number. Let $p_1, p_2, \cdots, p_l$ denote first $l$ primes and let $a_1, a_2, \cdots, a_l$ denote specific $l$ integers such that the **given** $x$ satisfies the following congruence system

$$x \equiv a_1 \mod(p_1)$$
$$x \equiv a_2 \mod(p_2)$$
$$\vdots$$
$$x \equiv a_l \mod(p_l)$$

We depict this fact in the form of sequence represented in the tabular form as follows:

| Mod | $p_1$ | $p_1$ | ... | $p_l$ |
|---|---|---|---|---|
| $x$ | $a_1$ | $a_2$ | ... | $a_l$ |

and call it the **remainder sequence** of $x$ with respect to primes $p_1, p_2, \cdots, p_l$.

**Algorithm for finding suitable $k$:**

Let $m < n$.
1. Create remainder sequences for $2m$ and $2n$ up to first $l$ primes $p_1(=2), p_2(=3), \cdots, p_l$ where $p_l^2 < 2n < p_{(l+1)}^2$ and let $2m \equiv \alpha_i \bmod(p_i)$ and $2n \equiv \beta_i \bmod(p_i)$.
2. Take $k \equiv \lambda_i \bmod(p_i)$ where $i = 1, 2, \cdots, l$ and $p_l^2 < 2n < p_{(l+1)}^2$ and determine properties that $k$ should obey by imposing conditions, namely, $\alpha_i \lambda_i$ as well as $\beta_i \lambda_i$ are **not** congruent to $1 \bmod(p_i)$.

**Remark 10.1:** Note that $k$ will be **appropriate** if $\alpha_i \lambda_i$ as well as $\beta_i \lambda_i$ are **not** congruent to $1 \bmod(p_i)$. Since, for example, if $\alpha_i \lambda_i \equiv 1 \bmod(p_i)$ then $2mk - 1$ will not be a prime. We take the smallest $k$ that satisfies the requirements present implicitly in the allowed values of remainders $\lambda_i$, $i = 1, 2, \cdots, l$ and extend the remainder sequences, if required, for $2m$, $2n$ and $k$ up to first $s$ primes $p_1(=2), p_2(=3), \cdots, p_s$ where $l \leq s$ and $p_s^2 < 2nk < p_{(s+1)}^2$.
3. Check whether this $k$ is appropriate. If YES stop. Else, if NO
4. Take next suitable $k$ and go to step 3. □

Note that among ($p_i$) possible values for $\lambda_i$ there will be some one value for which $\alpha_i \lambda_i \equiv 1 \bmod(p_i)$. Similarly, among ($p_i$) possible values for $\lambda_i$ there will be some one value for which $\beta_i \lambda_i \equiv 1 \bmod(p_i)$. Thus, there will be in all ($p_i - v_i$) values, (with $v_i = 0$, if $\alpha_i = 0$ and $\beta_i = 0$, and with $v_i = 1$ when $\alpha_i = \beta_i \neq 0$, so that $\alpha_i \lambda_i \equiv 1 \bmod(p_i)$ and $\beta_i \lambda_i \equiv 1 \bmod(p_i)$ holds for same $\lambda_i$, and $v_i = 2$ when $\alpha_i \neq \beta_i$ and both are nonzero) among the possible ($p_i$) values for $\lambda_i$. Thus, there are in all $\prod_{i=1}^{k}(p_i - v_i)$ possible values for $k$.

Now, using the Chinese remainder theorem we can generate all these numbers $k$ that are possible candidates so that the pair $(p, q) = (2mk - 1, 2nk - 1)$ will turn out to be a prime pair. In order to show the existence of a suitable candidate $k$ such that $p_l^2 < 2nk < p_{(l+1)}^2$ and $\alpha_i \lambda_i$ as well as $\beta_i \lambda_i$ are **not** congruent to $1 \mod(p_i)$ one could develop **similar arguments** as are given for the justification of Goldbach conjecture.

**Problem:** Develop **span argument** (as is done for the Goldbach conjecture) for showing the existence of a suitable $k$, to settle Schinzel-Sierpinski Conjecture.

**Example 10.1:** We express $\dfrac{m}{n} = \dfrac{11}{13} = \dfrac{(p+1)}{(q+1)}$ such that $p, q$ are prime numbers.

We construct **remainder sequence** for 22 and 26,

| Mod | 2 | 3 |
|---|---|---|
| 22 | 0 | 1 |

Here, primes up to $p_2$ are written in the first row and remainders modulo these primes for 22 are written in the second row below them.
Similarly,

| Mod | 2 | 3 |
|---|---|---|
| 26 | 0 | 2 |

Now we note down all the possibilities that exists for $\lambda_i$, $i = 1, 2$.
Clearly, $\lambda_1 = \{0,1\}$, $\lambda_2 = \{0\}$. Note that since $\lambda_2 = \{0\}$ therefore $k$ must be a multiple of three. When $k = 3$ we need to extend the remainder sequences up to fourth prime (since, $7 < \sqrt{26 \times 3} < 11$). With this the remainder sequences change to

| Mod | 2 | 3 | 5 | 7 |
|---|---|---|---|---|
| 22 | 0 | 2 | 2 | 1 |

| Mod | 2 | 3 | 5 | 7 |
|---|---|---|---|---|
| 26 | 0 | 2 | 1 | 5 |

So allowed values for $\lambda_3 = \{0,2,4\}$ and for $\lambda_4 = \{0,2,4,5,6\}$. But for $k = 3$ the remainder sequence is

| Mod | 2 | 3 | 5 | 7 |
|---|---|---|---|---|
| 3 | 1 | 0 | 3 | 3 |

i.e. when $k = 3$, $\lambda_3 = 3$ and $\lambda_4 = 3$, and these values are **disallowed** values. Also, when $k = 6$, $\lambda_3 = 1$ and $\lambda_4 = 6$, and among these the value $\lambda_3 = 1$ is a **disallowed** value. So, we next try k = 9. For this choice we need to extend the remainder sequences up to sixth prime, since $13 < \sqrt{26 \times 9} < 17$. Thus,

| mod | 2 | 3 | 5 | 7 | 11 | 13 |
|---|---|---|---|---|---|---|
| 22 | 0 | 2 | 2 | 1 | 0 | 9 |

and

| mod | 2 | 3 | 5 | 7 | 11 | 13 |
|---|---|---|---|---|---|---|
| 26 | 0 | 2 | 1 | 5 | 4 | 0 |

So allowed values for $\lambda_5 = \{0,2,4,5,6,7,8,9,10\}$ an for $\lambda_6 = \{0,1,2,4,5,6,7,8,9,10,11,12\}$. Also for k = 9 the remainder sequence becomes

| mod | 2 | 3 | 5 | 7 | 11 | 13 |
|---|---|---|---|---|---|---|
| 9 | 1 | 0 | 4 | 2 | 9 | 9 |

Thus, all values of $\lambda_i$ are **allowed values** in this case and so the primes we get are (197, 233) such that $\left(\dfrac{11}{13}\right) = \left(\dfrac{197+1}{233+1}\right)$.

**11. Mersenne, Fermat Primes and Other Twin Primes:** Consider the numbers of the form $M_q = 2^q - 1, q = 1,2,3,\cdots$. These numbers are called **Mersenne numbers** when the choice of $q$ is restricted to primes [5]. The following is the main conjecture related to Mersenne numbers: There are **infinitely many prime Mersenne numbers**.

If we fix a prime say $p > 2$ and consider congruence relations $2^q \equiv \alpha_q \mod(p)$ for all $q = 1,2,3,\cdots$ then it is easy to check the **important observation** that the remainders $\alpha_q$, $q = 1,2,3,\cdots$ form a **periodic sequence** of nonzero numbers with period $(p-1)$, and these numbers form some permutation of numbers $\{1,2,\cdots,p-1\}$. For $p = 2$ we get $\alpha_q = 0$ for all $q = 1,2,3,\cdots$

**Theorem 11.1:** The number of Mersenne primes $\leq x$ are **approximately** equal to $M(x) = \left(\dfrac{\log(x)}{\log(2)}\right) \prod_{p_j \leq \sqrt{x}} (1 - \dfrac{1}{(p_j - 1)})$.

**Proof:** A number of type $M_q \leq x$ is Mersenne prime if and only if $2^q \equiv \beta_i \mod(p_i)$ and $\beta_i \neq 1$ for all primes $2 < p_i \leq \sqrt{x}$. For a fixed $p_j > 2$ the remainder $\beta_j$ can take one of the values $\{1,\cdots,p_j - 1\}$. Therefore, for an $2^q \leq x$ chosen at random, the probability that $\beta_j \neq 1$ will be $\left(\dfrac{p_j - 2}{p_j - 1}\right)$. Clearly, for an $2^q \leq x$ **chosen at random**, the probability that $\beta_i \neq 1$ for all $i = 1,2,\cdots k$ such that $p_k^2 \leq x < p_{(k+1)}^2$ and $p_1, p_2, \cdots, p_k$ are all primes $\leq x$ will be

$\prod_{i=1}^{k}\left(\dfrac{p_i - 2}{p_i - 1}\right) = \prod_{i=1}^{k}\left(1 - \dfrac{1}{(p_i - 1)}\right)$. Now, the numbers of type $2^q \leq x$ equal to $u$ where $u = \max\{q\}$ such that $2^u \leq x$. It is easy to check that $u = \left(\dfrac{\log(x)}{\log(2)}\right)$. Therefore, the number of Mersenne primes less than or equal to $x$ will be

$\left(\dfrac{\log(x)}{\log(2)}\right)\prod_{i=1}^{k}\left(1 - \dfrac{1}{(p_i - 1)}\right) = \left(\dfrac{\log(x)}{\log(2)}\right)\prod_{p_i \leq \sqrt{x}}\left(1 - \dfrac{1}{(p_i - 1)}\right) = M(x) \square$

It is clear to see that $M(x) \to \infty$ as $x \to \infty$. Thus, Mersenne primes are **infinite!**

**Remark 11.1:** One can develop a **combinatorial formula** to count the Mersenne primes less than or equal to $x$ as is done for counting primes less than or equal to $x$ in lemma 2.1.

**Theorem 11.2:** Let $S = \{p_1, p_2, \cdots, p_k\}$ be set of first $k$ primes and $p_k^2 \leq x < p_{(k+1)}^2$, then the number of Mersenne primes less than or equal to $x$ can be given by the following formula

$$M(x) = u - \sum_{i=2}^{k}\left[\frac{u}{p_i(1)}\right] + \sum_{i<j}^{k}\left[\frac{u}{p_i(1)p_j(1)}\right] - \sum_{i<j<l}^{k}\left[\frac{u}{p_i(1)p_j(1)p_l(1)}\right] + \cdots$$

$$\cdots + (-1)^k\left[\frac{u}{p_1(1)p_2(1)\cdots p_k(1)}\right] + (\lambda_k - 1)$$

where [q] represents the integral part of q, $u = \left[\frac{\log(x)}{\log(2)}\right]$ and $\left[\frac{u}{p_j(1)}\right]$ represents the numbers of type $2^q \leq x$ such that $2^q \equiv 1 \bmod(p_j)$ and where $\lambda_k$ stands for the Mersenne primes $\leq \sqrt{x}$.

**Proof:** Follows using inclusion-exclusion principle. □

We now proceed with the discussion of the other number called the Fermat numbers. This time, consider the numbers of the form $F_q = 2^q + 1, q = 1, 2, 3, \cdots$. These numbers are called **Fermat numbers** when the choice of $q$ is restricted to the numbers of the form $2^n, n \geq 0$, [5]. The following is the main conjecture related to Fermat numbers:

There are **infinitely many prime Fermat numbers**. This result has a special significance because of the crowning achievement of Gauss: If a regular polygon may be constructed using ruler and compass having $n(\geq 3)$ sides then $n = 2^k p_1 p_2 \cdots p_l$, where $k, l \geq 0$, and $p_1, p_2, \cdots, p_l$ are distinct odd Fermat primes, (see page 74, [5]).

**Theorem 11.3:** The number of Fermat primes $\leq x$ are **approximately** equal to $F(x) = \left(\dfrac{\log(x)}{\log(2)}\right) \prod_{p_j \leq \sqrt{x}} (1 - \dfrac{1}{(p_j - 1)})$.

**Proof:** A number of type $F_q \leq x$ is Fermat prime if and only if $2^q \equiv \beta_i \mod(p_i)$ and $\beta_i \neq (p_i - 1)$ for all primes $2 < p_i \leq \sqrt{x}$. For a fixed $p_j > 2$ the remainder $\beta_j$ can take one of the values $\{1, \cdots, p_j - 1\}$. Therefore, for an $2^q \leq x$ chosen at random, the probability that $\beta_j \neq (p_j - 1)$ will be $\left(\dfrac{p_j - 2}{p_j - 1}\right)$. Clearly, for an $2^q \leq x$ **chosen at random**, the probability that $\beta_i \neq (p_i - 1)$ for all $i = 1, 2, \cdots k$ such that $p_k^2 \leq x < p_{(k+1)}^2$ and $p_1, p_2, \cdots, p_k$ are all primes $\leq x$ will be $\prod_{i=1}^{k} \left(\dfrac{p_i - 2}{p_i - 1}\right) = \prod_{i=1}^{k} \left(1 - \dfrac{1}{(p_i - 1)}\right)$. Now, the numbers of type $2^q \leq x$ equal to $u$ where $u = \max\{q\}$ such that $2^u \leq x$. It is easy to check that $u = \left(\dfrac{\log(x)}{\log(2)}\right)$. Therefore, the number of Fermat primes less than or equal to $x$ will be

$$\left(\dfrac{\log(x)}{\log(2)}\right) \prod_{i=1}^{k}\left(1 - \dfrac{1}{(p_i - 1)}\right) = \left(\dfrac{\log(x)}{\log(2)}\right) \prod_{p_i \leq \sqrt{x}} \left(1 - \dfrac{1}{(p_i - 1)}\right) = F(x) \square$$

It is clear to see that $F(x) \to \infty$ as $x \to \infty$. Thus, Fermat primes are **infinite!**

**Remark 11.2:** One can develop a **combinatorial formula** to count the Fermat primes less than or equal to $x$ as is done for counting primes less than or equal to $x$ in lemma 2.1.

**Theorem 11.4:** Let $S = \{p_1, p_2, \cdots, p_k\}$ be set of first $k$ primes and $p_k^2 \leq x < p_{(k+1)}^2$, then the number of Fermat primes less than or equal to $x$ can be given by the following formula

$$M(x) = u - \sum_{i=2}^{k}\left[\frac{u}{p_i(\chi)}\right] + \sum_{i<j}^{k}\left[\frac{u}{p_i(\chi)p_j(\chi)}\right] - \sum_{i<j<l}^{k}\left[\frac{u}{p_i(\chi)p_j(\chi)p_l(\chi)}\right] + \cdots$$

$$\cdots + (-1)^k\left[\frac{u}{p_1(\chi)p_2(\chi)\cdots p_k(\chi)}\right] + (\lambda_k - 1)$$

where [q] represents the integral part of q, , $u = \left[\dfrac{\log(x)}{\log(2)}\right]$ and

$\left[\dfrac{u}{p_j(\chi)}\right]$ represents the numbers of type $2^q \leq x$ such that

$2^q \equiv (\chi)\mod(p_j)$ and $\chi = (p_j - 1)$. Note that $\lambda_k$ stands for the Fermat primes $\leq \sqrt{x}$.

**Proof:** Follows using inclusion-exclusion principle. □

**Remark 11.3:** The important observation stated above about the periodicity of the remainders $\alpha_q$, $q = 1,2,3,\cdots$ while one considers the congruence relations $2^q \equiv \alpha_q \mod(p)$ for all $q = 1,2,3,\cdots$ hold good even for the congruence relations $(k \times a^q) \equiv \alpha_q \mod(p)$ for all $q = 1,2,3,\cdots$, where $k, a$ are arbitrary positive integers. So, one can easily generalize the above theorems about Mersenne and Fermat numbers to the numbers of the form $k \times a^q \pm r$, where $r$ is some known fixed constant, e.g. $r = 1, 3, 5, \ldots$ etc.

**Remark 11.4:** One more important conjecture related to Mersenne numbers is as follows: There are **infinitely many composite Mersenne numbers**. The settlement of this conjecture requires infinitude of a special kind of Sophie Germain primes $q$ of the form $q = k \times 2^n - 1$ (so, $2q + 1$ is also prime). The existence of such a pair yields a composite Mersenne number $M_q = 2^q - 1$ by the classical result stated by Euler and proved by Lagrange and again by Lucas, (see page 76 very admirable book, [5]). As stated in the above remark 11.3 one can achieve the desired infinitude of these special kind of Sophie Germain primes by proceeding on the similar lines as is done for theorems 11.1, 11.2.

**Other Twin Primes:** Prime numbers are infinite since the time when Euclid gave his one of the most beautiful proof of this fact! Prime number theorem (PNT) reestablishes this fact and further it also gives estimate about the count of primes less than or equal to $x$. PNT states that as $x$ tends to infinity the count of primes up to $x$ tends to $x$ divided by the natural logarithm of $x$. Twin primes are those primes $p$ for which $p+2$ is also a prime number. The well known twin prime conjecture (TPC) states that twin primes are (also) infinite. Related to twin primes further conjectures that can be made by extending the thought along the line of TPC, are as follows: Prime numbers $p$ for which $p+2n$ is also prime are (also) infinite for all $n$, where $n = 1(\text{TPC}), 2, 3, \ldots, k, \ldots$. Now, we provide a simple argument in support of all twin prime conjectures.

The celebrated prime number theorem (PNT) gives exact estimate for cardinality of primes up to $x$. If $\pi(x)$ denotes the number of primes less than or equal to $x$ then

$$\pi(x) = \frac{x}{\ln(x)} \quad \text{as } x \to \infty.$$

Let us denote by $\pi_{2n}(x), n = 1, 2, \cdots, k, \cdots$ the number of primes, $p$, less than or equal to $x$ for which $p+2n$ is also a prime number. $\pi_{2n}(x)$ may be appropriately called as twin primes of $n$-type. Clearly, twin primes of 1-type are usual twin prime numbers.

**All kinds of Twin Primes are Infinite:** According to prime number theorem (PNT) when $x$ is tending to infinity the cardinality of primes up to $x$, $\pi(x)$, is given by $\pi(x) = \frac{x}{\ln(x)}$.

Let us consider all possible distinct prime pairs made out of all prime numbers up to $x$. So, let $\{p_1, p_2, \cdots, p_k\}$ be the primes up to $x$. Then initially we will form following types of pairs of primes $\{(p_1, p_1), (p_2, p_2), \cdots, (p_k, p_k)\}$. The count of these pairs will be obviously equal to $\pi(x)$. We then will form the pairs of primes $\{(p_1, p_2), (p_2, p_3), \cdots, (p_{k-1}, p_k)\}$, then in continuity we will form the pairs of primes $\{(p_1, p_3), (p_2, p_4), \cdots, (p_{k-2}, p_k)\}, \ldots$, finally we form $\{(p_1, p_k)\}$. Obviously, if we will consider the count of all these pairs of distinct prime together then it will be clearly equal to $\binom{\pi(x)}{2}$. Thus, the total count of all prime pairs will be will be

$$\pi(x) + \binom{\pi(x)}{2} = \frac{1}{2}\pi(x)(1+\pi(x)) = \frac{1}{2}[\pi(x) + (\pi(x))^2] \quad \ldots(A)$$

It is very interesting to observe that actually the count of these pairs is also equal to

$$\pi(x) + \pi_2(x) + \pi_4(x) + \cdots \pi_{2n}(x) + \cdots \quad \ldots(B)$$

thus, equating these two counts, (A) and (B), and dividing by both sides by $\pi(x)$ we
have

$$\frac{1}{2}[1 + \pi(x)] = [1 + \frac{\sum_{n=1}^{\infty} \pi_{2n}(x)}{\pi(x)}] \quad \ldots(C)$$

as $x \to \infty$.

**Example:** Let us consider primes up to first 50 integers.

**Primes:** {3, 5, 7, 11, 13, 17, 19, 23, 29, 31, 37, 41, 43, 47}.

Therefore, $\pi(50) = 14$.

So, $\frac{1}{2}[\pi(50)][1 + \pi(50)] = 105$.

Also, .

$$[\pi(50) + \pi_2(50) + \pi_4(50) + \pi_6(50) + \cdots + \pi_{44}(50)] = 105$$
.

As it is equal to:

$$\begin{bmatrix} 14 + 6 + 6 + 9 + 5 + 6 + 7 + 5 + 4 + 6 + 4 + 2 + \\ 6 + 4 + 3 + 4 + 2 + 3 + 3 + 2 + 2 + 1 + 1 \end{bmatrix} = 105$$.

It is clear to see that the LHS of above equation (C) goes to infinity $x \to \infty$, since as per PNT $\pi(x) = \dfrac{x}{\ln(x)}$ which clearly goes to infinity as $x \to \infty$. On RHS of equation (C) we have collection of terms made up of densities of twin primes of various types, $\pi_{2n}(x), n = 1, 2, \cdots, k, \cdots$ divided by density of primes, $\pi(x)$. So, if the count of each term on RHS of equation (C), $\pi_{2n}(x), n = 1, 2, \cdots, k, \cdots$, is finite and converges to some finite value as $x \to \infty$ then RHS of equation (C) will be made up of terms which are all equal to zero (except the first term which is equal to 1). Therefore RHS of equation (C) can't diverge to infinity. Actually, RHS of equation (C) should diverge to infinity to avoid contradiction since LHS of equation (C) contains $\pi(x)$ and so is diverging to infinity when $x \to \infty$. Now on RHS of equation (C) since we have a divergent quantity, namely, the density of total primes, $\pi(x)$, present in the denominator as mentioned above, therefore, each $\pi_{2n}(x), n = 1, 2, \cdots, k, \cdots$ must be diverging to infinity, or some infinitely many of $\pi_{2n}(x), n = 1, 2, \cdots, k, \cdots$, for some infinitely many $n$, must be diverging to infinity. But since count of primes counted in the each of the twin prime densities $\pi_{2n}(x), n = 1, 2, \cdots, k, \cdots$ is actually a subset of $\pi(x)$ therefore each, or, at least some infinitely many of $\pi_{2n}(x), n = 1, 2, \cdots, k, \cdots$, though should diverge to infinity still must be diverging slowly than $\pi(x)$ and therefore as $x \to \infty$ each ratio $\dfrac{\pi_{2n}(x)}{\pi(x)}, n = 1, 2, \cdots, k, \cdots$ must be converging to a (nonzero) positive number, $\lambda_n$, such that series $1 + \sum_{n=1}^{\infty} \lambda_n$ should form a divergent series as required for maintaining consistency with the divergent LHS of equation (C)! Thus, at least some infinitely many of $\pi_{2n}(x), n = 1, 2, \cdots, k, \cdots$, though slowly than $\pi(x)$ must be diverging to infinity as $x \to \infty$ !!

## Acknowledgements

The author is thankful to Dr. M. R. Modak, Bhaskaracharya Pratishthana, Pune, for useful discussions.## References

1. T. L. Heath, The Thirteen Books of Euclid's Elements, Book IX, Proposition 20, Dover Publications, New York, 1956.
2. V. K. Balakrishnan, Schaum's Outline Series, Theory and Problems of Combinatorics, McGraw-Hill, Inc, 1995.
3. Tom M. Apostol, Introduction to Analytic Number Theory, Springer International Student Edition, Narosa Publishing House, New Delhi, 1980.
4. Ivan Niven, Herbert S. Zuckerman, Hugh L. Montgomery, An Introduction to The Theory of Numbers, Fifth Edition, John Wiley & Sons, Inc, 2004.
5. Paulo Ribenboim, The Little Book of Bigger Primes, Springer International Edition (Second Edition), Springer (India) Private Limited, New Delhi─110 001, India, 2005.
6. M. L. Perez, Generalizations of Goldbach and Polignac conjectures, www.gallup.unm.edu/~smarandache/prim-sum.txt
7. Melvyn B. Nathanson, Methods in Number Theory, Page 288, Open Problem 5, Springer- Verlag, New York, Inc, First Indian Reprint, 2005.50